\documentclass[12pt,a4paper]{amsart}

\usepackage[english]{babel}
\usepackage{amssymb,graphicx,bbm,xcolor,calc,ifthen,caption}
\usepackage[a4paper, total={6in, 8in}]{geometry}

\catcode`\<=\active \def<{
\fontencoding{T1}\selectfont\symbol{60}\fontencoding{\encodingdefault}}
\catcode`\>=\active \def>{
\fontencoding{T1}\selectfont\symbol{62}\fontencoding{\encodingdefault}}
\catcode`\|=\active \def|{
\fontencoding{T1}\selectfont\symbol{124}\fontencoding{\encodingdefault}}
\newcommand{\assign}{:=}
\newcommand{\dueto}[1]{\textup{\textbf{(#1) }}}
\newcommand{\mathd}{\mathrm{d}}
\newcommand{\nocomma}{}
\newcommand{\tmem}[1]{{\em #1\/}}
\newcommand{\tmname}[1]{\textsc{#1}}
\newcommand{\tmop}[1]{\mathrm{#1}\hspace{.5em}}
\newcommand{\tmstrong}[1]{\textbf{#1}}
\newcommand{\tmtextit}[1]{{\itshape{#1}}}
\newenvironment{proof*}[1]{\noindent\textbf{#1\ }}{\hspace*{\fill}$\Box$\medskip}
\newtheorem{proposition}{Proposition}
\newtheorem{theorem}{Theorem}
\newcommand{\tmfloatcontents}{}
\newlength{\tmfloatwidth}
\newcommand{\tmfloat}[5]{
  \renewcommand{\tmfloatcontents}{#4}
  \setlength{\tmfloatwidth}{\widthof{\tmfloatcontents}+1in}
  \ifthenelse{\equal{#2}{small}}
    {\ifthenelse{\lengthtest{\tmfloatwidth > \linewidth}}
      {\setlength{\tmfloatwidth}{\linewidth}}{}}
    {\setlength{\tmfloatwidth}{\linewidth}}
  \begin{minipage}[#1]{\tmfloatwidth}
    \begin{center}
      \tmfloatcontents
      \captionof{#3}{#5}
    \end{center}
  \end{minipage}}
\newcommand{\cal}[1]{\mathcal{#1}}

\begin{document}

\title[Existence of spray flames travelling waves ]{Existence of travelling waves and high activation energy limits for a
onedimensional thermo-diffusive lean spray flame model}

\author{Pierre Berthonnaud$^1$}
\address{$^1$Institut de Math\'ematiques de Toulouse. 118, route de Narbonne. 31062 Toulouse}

\author{Komla Domelevo$^1$}

\begin{abstract}
  We provide a mathematical analysis of a thermo-diffusive combustion model
  of lean spray flames, for which we prove the existence of travelling waves.
  In the high activation energy singular limit we show the existence of two
  distinct combustion regimes with a sharp transition -- the diffusion limited
  regime and the vaporisation controlled regime. The latter is specific to
  spray flames with slow enough vaporisation. We give a complete
  characterisation of these regimes, including explicit velocities, profiles,
  and upper estimate of the size of the internal combustion layer.
  
  Our model is on the one hand simple enough to allow for explicit asymptotic
  limits and on the other hand rich enough to capture some particular aspects
  of spray combustion. Finally, we briefly discuss the cases where the
  vaporisation is infinitely fast, or where the spray is polydisperse.
\end{abstract}

\maketitle

\vspace{2pc}
\noindent{\it Keywords}: reaction-diffusion equations, spray flames, travelling
waves, high activation energy, singular limits

\section{Introduction}

This paper provides a rigourous mathematical analysis of some aspects of spray
combustion, including the analysis of the so-called high activation energy
limit for a spray flame model. This notion of high activation energy limit was
first introduced in the pioneering work of Zeldovich and Frank-Kamenetskii
{\cite{ZelFra1938}}, and refers to the limit where the combustion rate is much
faster than any other physical phenomenon, in particular diffusion. Since then, there
have been many studies of these asymptotics and applications to gaseous flames.
However, gas-vapor-droplets systems have many applications in industry or everyday's
life, such as diesel or propulsion engines. When trying to understand some of their
specific features, it appears that the structure of these two-phase flames as
well as their speed or stability are greatly affected by the presence of
vaporising liquid droplets possibly interacting with the combustion zone.

The behaviour of spray flames has been investigated a lot in the literature
and a wide variety of regimes were considered. Dating from the 70s and early
80s, we can quote the works of {\tmname{Polymeropoulos}} et al.
{\cite{Pol1974a}}{\cite{Pol1975a}}, {\tmname{Mizutani}} et al.
{\cite{MizNak1973a}}{\cite{MizNak1973b}}, {\tmname{Hayashi}} et al.
{\cite{HayKum1974a}}{\cite{HayKumSak1976a}}. {\tmname{Ballal}} et al.
{\cite{BalLef1981a}}. They present studies of the propagation of liquid fuel
sprays, study the influence of the size of the droplets, the type and geometry
of the spray flame as well as its structure. The investigation of those
vapor--drop--air systems was continued in the 90s. To quote but a few works,
we refer to the work of \tmname{Aggarwal} and \tmname{Sirignano} {\cite{AggSir1985a}} as well as
the papers of {\tmname{Greenberg}}, {\tmname{Tambour}} and
{\tmname{Silverman}} {\cite{SilGreTam1991}}{\cite{SilGreTam1992a}}, where the
structure of spray flames as well as the influence of parameters such as
droplet size, fuel volatility, or equivalence ratio are also investigated
analytically. More elaborate situations appearing in propulsion engines are
for example pulsating or acoustic instabilities. We refer to
{\cite{GreMcIBri1999}}{\cite{GreMcIBri2001}}{\cite{Gre2002a}}{\cite{Gre2003}}{\cite{GreDvo2003}}{\cite{ClaSun1990}}
and more recently
{\cite{ConDomRoqRyz2006a}}{\cite{GreKagSiv2009a}}{\cite{KagGreSiv2010a}}{\cite{KagGreSiv2012a}}{\cite{KatGre2013a}}
for studies in that direction.

Only few of these studies of spray flames involve
rigorous--in--the--mathematical--sense analysis of spray flames models where
existence, uniqueness, or asymptotic limits are derived (see e.g.
{\cite{LauMas2002b,ConDomRoqRyz2006a}}). This is in contrast with purely
gaseous combustion, where a lot of results exist in the mathematical
literature for various regimes and asymptotics. In particular, a complete
study of thermo-diffusive lean gaseous flame fronts and the high activation
energy limit is in the paper of \tmname{Berestycki}, \tmname{Nicolaenko} and \tmname{Scheurer}
{\cite{BerNicSch1985}}, and in the paper of \tmname{Berestycki} and \tmname{Larrouturou}
{\cite{BerLar1991}}. The originality of the present work is to provide a
complete mathematical analysis of the counterpart of those systems, namely a
lean {\tmstrong{spray}} flame model that on the hand is {\tmstrong{simple
enough}} to allow for explicit asymptotic limits, and on the other hand
{\tmstrong{rich enough}} to capture some particular aspects of spray
combustion. Those results should be compared with the very interesting work of
\tmname{Suard}, \tmname{Nicoli} and \tmname{Haldenwang} {\cite{SuaNicHal2001}}. Using numerical
simulations, these authors investigate the scaling laws of the spray flame
with respect to the vaporisation rate of the liquid phase. When the
vaporisation is fast enough, the velocity of the spray flame is comparable to
that of the gaseous flame where the whole reactant would be present only in a
gaseous form. This is the {\tmstrong{diffusion controlled}} regime. On the
contrary when the time for complete vaporisation exceeds a critical value, the velocity
of the spray flame starts to decrease. This is the so--called
{\tmstrong{vaporisation controlled}} regime. In this case, they also investigated
the structure of the reaction zone, and showed that it has a more complex
structure than a gaseous flame.

In this work, we prove the existence of spray flame travelling waves and the
existence of {\tmstrong{two distinct combustion regimes}} with a sharp
transition in the {\tmstrong{high activation energy}} ({\tmstrong{HAE}}). In
this limit, we provide a complete characterisation of the profiles, that can
be written explicitly provided simple analytic expressions of the vaporisation law. We also
show that the internal combustion layer is likely to be much larger in the
vaporisation controlled regime as compared to the size of the internal
combustion layer of a comparable gaseous flame.

\subsection{The equations of the lean spray flame model}

\subsubsection*{The classical gaseous reactive flow.}A gas mixture of total mass
density $\rho$ is considered to be made of $N$ fluids corresponding to the
different species present in the mixture. If the mass density of those species
are $( \rho_{i} )_{i=1,N}$ then by definition $\rho = \sum_{i} \rho_{i}$, and
the respective mass fractions $( Y_{i} )_{i=1,N}$ of the species are defined
by $\rho_{i} = \rho Y_{i}$. In the case of a single reactant $Y$, the chemical
reaction writes as
\[ A \hspace{1em} \rightarrow \hspace{1em} P, \]
where $A$ denotes the reactant species and $P$ denotes one (or a linear
combination) of the species produced by the reaction. Notice that the rate of
variation of the single reactant $A$ is directly proportional to minus that of
the products and it suffices to determine the mass fraction $Y$ of the
reactant $A$ in order to determine the mass fractions of the other species.

The mass, momentum, energy equation for a reactive flow in one dimension
without external forces write, setting $D_{t} \assign \partial_{t} +v
\partial_{x}$,
\[ \partial_{t} \rho + \partial_{x} ( \rho v ) =D_{t} \rho + \rho \partial_{x}
   v=0, \]
\[ \rho D_{t} v+ \partial_{x} p= \frac{4}{3} \partial_{x} ( \kappa
   \partial_{x} v ) , \]
\[ \rho c_{p} D_{t} T- \partial_{x} ( \lambda \partial_{x} T ) =Q \omega
   +D_{t} p+ \frac{4}{3} \kappa ( \partial_{x} v )^{2} , \]
where $v$ is the mass--average velocity of the mixture, $p$ the hydrostatic
pressure, $\kappa$ the dynamic viscosity coefficient, $T$ the temperature,
$c_{p}$ the specific heat at constant pressure, $\lambda \assign \lambda ( T
)$ the thermal conductivity, $Q$ the chemical heat release of the reaction and
$\omega$ is the rate at which the reaction occurs and has the form
\[ \omega =B ( T ) \frac{\rho Y}{\mu} \exp \left( - \frac{E}{R T} \right) , \]
where $E$ is the activation energy of the single reaction, $R$ the perfect gas
constant, $T_{A} =E/R$ the activation energy temperature, $\mu$ the molecular
mass of the reactant $A$, and $B ( T )$ a non--stiff prefactor. Finally the
equation for the mass fraction $Y$ of the reactant $A$ obeys
\[ \rho D_{t} Y- \partial_{x} ( \gamma \partial_{x} Y ) =- \mu \omega . \]
We refer the reader to {\cite{Smo1983,Wil1988}} for more details.

\subsubsection*{The liquid phase.} In the situation where the reactant $A$ is also
present in the mixture in the form of vaporising liquid droplets that are well
dispersed in the mixture, one can treat the liquid phase as a new continuous
``species'', where the droplets are homogeneously spread inside the mixture,
leading to both homogeneized vaporisation and combustion in the bulk of the
gaseous phase. This is certainly not satisfied in practice, where spray
combustion can typically involve flames surrounding individual droplets, with
strong temperature gradients within the flame. However, we need this
hypothesis for the derivation of our model. Also, ideally the continuity
equation for the gaseous mixture should take into account the volume fraction
occupied by the liquid phase, the momentum equation should incorporate the
drag forces, the energy equation should include the loss of heat due to the
latent heat of the liquid as well as the dissipation due to the drag forces,
etc.

We are going to consider regimes in which we can discard those couplings.
Since the mass density of the liquid reactant $A$ is in general several orders
of magnitude larger than that of the gaseous species, on can assume that the volume fraction of
the liquid phase is small. The drag forces are always present and can produce
complex turbulent flows. However, the drag force is typically inversely
proportional to the surface area of the droplets and proportional to the
difference of velocities of the droplets and the surrounding gas. Small enough
droplets behave as ``passive scalars'', that is particles whose velocity can
be set equal to that of the surrounding fluid. The drag forces are internal
forces for the two-phase flow and do not affect the total momentum $( \rho_{l}
v_{l} + \rho v )$ of the liquid--gas mixture. As a consequence, since $v_{l}
\simeq v$, one can neglect the drag forces in the gaseous mixture if the mass
density $\rho_{l}$ of the liquid phase is small enough compared to the mass
density $\rho$ of the gas mixture. The same conclusion holds for the energy
equation.

Let now $M$ denote the mass of an isolated vaporising droplet immersed in a
gas. The vaporisation law of the droplet can be very complex
(see e.g. {\cite{Wil1958a,Wil1976,ORo1981,Wil1988,Sir1999}}). As an
approximation, we will assume that the vaporisation rate $\phi$ only depends
upon the temperature $T$ of the surrounding gas and the mass $M$ of the
droplet , i.e.
\[ \mathd_{t} M=- \phi ( T,M ) . \]
We will be mainly interested in the monodisperse case, where all droplets in
the unburnt gas have the same size. In the laminar flows with constant velocity that we will be considering,
this implies that particles located at the same position $x$ at time $t$ all
have the same size. We can therefore introduce the mass profile $M ( t,x )$
which obeys,
\[ D_{t} M=- \phi ( T,M ) . \]
The situation for polydisperse sprays is briefly considered in
the HAE limit in Section \ref{S: extensions}. Let further $n ( t,x )$ be
the number density of droplets, that is the number of droplets per unit
volume. The mass density of the liquid phase writes $\rho_{l} ( t,x ) =n ( t,x
) M ( t,x )$. In the laminar flows we are considering, the droplets do not
coalesce nor break--up, hence
\[ D_{t} n=0, \]
that is the number density of the droplets is simply advected by the flow, and
$D_{t} \rho_{l} =n ( t,x ) D_{t} m$. It follows that the continuity equation
for the gaseous reactant $A$ now incorporates the vaporisation flux as a
source term
\[ \rho D_{t} Y- \partial_{x} ( \gamma \partial_{x} Y ) =- \mu \omega -D_{t}
   \rho_{l} =- \mu \omega -n D_{t} M. \]

Finally, we assume that the effects of the latent heat are negligible. This is
true when the mass density of the liquid phase is small enough compared to
that of the gaseous mixture, which we assumed also in order to neglect the
drag forces. The study of the influence of the latent heat is certainly an
important and interesting feature since it can greatly influence the speed or
existence of the flame. However this study is beyond the scope of the present
work.

\subsubsection*{Two--phase system.} We are interested in laminar, low Mach number,
regimes where one can neglect the pressure terms as well as the kinematic
viscosity. We want to study the existence and properties of travelling waves
moving at constant speed $v_{0}$ to the left. In the frame of reference of an
observer moving at speed $v_{0}$ to the left, the profiles are steady,
solution to the system
\[ v_{0} \rho' + ( \rho v )' =0, \]
\[ \rho v_{0} v' =0, \]
\[ c_{p} ( \rho v_{0} + \rho v ) T' - ( \lambda T' )' =Q \omega , \]
\[ n' =0, \]
\[ ( \rho v_{0} + \rho v ) Y' - ( \gamma Y' )' =- \mu \omega -n  ( v_{0} +v )
   M' . \]
From the continuity equation, the mass flux $c \assign \rho ( v_{0} +v )$ is a
constant. Moreover $v' =0$ implies $v=0$ in the moving frame, and therefore
also $\rho' =0$, so that the gas density is constant, $\rho \assign \rho_{0}$.
Setting also $n:=n_0$, a constant, the
system reduces to
\[ c c_{p} T' - ( \lambda T' )' =Q \omega , \]
\[ c Y' - ( \gamma Y' )' =- \mu \omega -c  \frac{n_{0}}{\rho_{0}}  M' , \]
\[ c  \frac{n_{0}}{\rho_{0}}  M' =-n_{0}   \phi ( T,M ) , \]
with prescribed values $(T_u,Y_u,n_u,M_u)$ in the unburnt region,
\[ T ( - \infty ) =T_{u} , \hspace{1em} Y ( - \infty ) =Y_{u} , \hspace{1em}
   n_{0} =n_{u} , \hspace{1em} M ( - \infty ) =M_{u} , \]
and assuming complete reaction and complete vaporisation on the burnt region
\[ Y ( + \infty ) =0, \hspace{1em} M ( + \infty ) =0 \nocomma . \]
The temperature $T ( + \infty ) \assign T_{b}$ in the burnt region is obtained
by integrating the equation for the quantity
\[ \frac{c_{p}}{Q}  T+ \frac{1}{\mu}   \left( Y+ \frac{n_{0}}{\rho_{0}}  M
   \right) \]
on $( - \infty ,+ \infty )$. It follows
\[ T_{b} = \frac{Q}{c_{p} \mu} \left( Y_{u} + \frac{n_{0}}{\rho_{0}}  M_{u}
   \right) . \]

\subsubsection*{Normalized system.}The normalized variables are $( u,v,m )$ defined as
\begin{equation}
  u= \frac{T-T_{u}}{T_{b} -T} , \hspace{1em} v= \frac{\rho_{0} Y}{\rho_{0}
  Y_{u} +n_{0} M_{u}} \nocomma , \hspace{1em} m= \frac{M}{\rho_{0} Y_{u}
  +n_{0} M_{u}} , \label{eq: normalized variables}
\end{equation}
solution to the system,
\begin{equation}
  \begin{array}{rll}
    -u'' +c u' & = & \tilde{f} ( u ) v \hspace{1em} \tmop{on}  \mathbbm{R},\\
    - \Lambda v'' +c v' & = & - \tilde{f} ( u ) v-c n_{0} m'  _{} \hspace{1em}
    \tmop{on}  \mathbbm{R},\\
    c m' & = & - \tilde{\phi} ( u,m ) \hspace{1em} \tmop{on}  \mathbbm{R},\\
    u ( - \infty ) =0, &  & u ( + \infty ) =1,\\
    v ( - \infty ) =v_{u} , &  & v ( + \infty ) =0,\\
    m ( - \infty ) =m_{u} , &  & 
  \end{array} \label{eq: nondimensional travelling wave}
\end{equation}
where $\Lambda = \gamma c_{p} / \lambda$ is the reciprocal of the Lewis
number, and with appropriate renormalized formulas for the reaction rate
$\tilde{f}$ and the vaporisation rate $\tilde{\phi}$ (we omit the tildes
below). Also $n_{0} \geqslant 0$, $v_{u} \geqslant 0$, $m_{u} \geqslant 0$ and
$v_{u} +n_{0} m_{u} =1$.

\subsubsection*{Reaction and vaporisation laws.}In order to avoid the cold boundary
effect, we assume that the reaction rate is zero below the normalized ignition
temperature $0< \theta_{i} <1$, namely $f: [ 0,1 ] \rightarrow
\mathbbm{R}^{+}$, $f ( u ) =0$, for all $0 \leqslant u< \theta_{i}$, $f$
positive on $( \theta_{i} ,1 ]$ and Lipschitz continuous on $[ \theta_{i} ,1
]$.

For the study of the high activation energy limit, we will use on $(
\theta_{i} ,1 ]$ the normalized arrhenius law
\begin{equation}
  f_{\varepsilon} ( u ) \assign \frac{1}{\varepsilon^{2}}   \exp \left(
  \frac{u-1}{\varepsilon} \right) , \label{eq: Arrehnius law}
\end{equation}
where $\varepsilon$ denotes the inverse of the (normalized) activation
temperature. We also introduce $\mu$, $0<\mu<+\infty$, defined as

\[ \mu \assign \lim_{\varepsilon \rightarrow 0} \int_{\theta_{i}}^{1}
   f_{\varepsilon} ( s ) ( 1-s ) \mathd s. \]

We assume that vaporisation starts at temperature $\theta_{v} >0$, i.e. $\phi
( u,m ) =0$ for all $0 \leqslant u< \theta_{v} \nocomma$ and we also impose
$\phi ( u,0 ) =0$ for all $0 \leqslant u \leqslant 1$. It is also natural to
assume that $\theta_{v} < \theta_{i}$ and to impose $\phi : [ 0,1 ] \times [
0,1 ] \rightarrow \mathbbm{R}^{+}$ is an increasing function of $u$, a
decreasing function of $m$, and is a positive Lipschitz function on $[
\theta_{v} ,1 ] \times [ 0,1 ]$. Let be given any temperature $\theta$, with $1>
\theta > \theta_{v}$, and any mass $m$. We finally make the natural assumption
that the interval of time $\tau ( \theta ,m )$ needed for complete
vaporisation of a droplet of initial mass $m$ at temperature $\theta$ is
finite.

\subsection{Main results and summary}Section \ref{S: Existence of travelling
waves} is devoted to the proof of the existence of travelling waves solutions
to system (\ref{eq: nondimensional travelling wave}). Precisely,

\begin{theorem}
  \label{T: existence of travelling waves}{\dueto{existence of travelling
  waves}}System (\ref{eq: nondimensional travelling wave}) admits a solution
  in $X= \mathcal{C}^{1} ( \mathbbm{R} ) \times \mathcal{C}^{1} ( \mathbbm{R}
  ) \times \mathcal{C}^{0} ( \mathbbm{R} ) \times \mathbbm{R}$.
\end{theorem}

This result involves \tmtextit{a priori} bounds that are proved in Section
\ref{SS: preliminary estimates}, and where we emphasize the technical
difficulties induced by the presence of vaporising droplets. This allows us to
prove in Section \ref{SS: existence of a solution on bounded domains} the
existence of solutions to a system similar to (\ref{eq: nondimensional
travelling wave}) but defined on a bounded domain with appropriate boundary
conditions. The proof of Theorem \ref{T: existence of travelling waves}
follows in Section \ref{SS: Existence of travelling waves}.

Section \ref{S: High activation energy limit} is devoted to the asymptotic
analysis of the system in the high activation energy limit. In Section
\ref{SS: characterisation of the limiting profiles}, we prove the convergence
of the system towards a Dirac model for spray flames:

\begin{theorem}
  \label{T: limiting system in the HEA limit}{\dueto{the Dirac model for spray
  flames}}Let the reaction term $f_{\varepsilon} ( \cdot )$ as in (\ref{eq:
  Arrehnius law}). There exists a decreasing sequence $( \varepsilon_{n} )_{n
  \in \mathbbm{N}}$ such that $( u_{\varepsilon_{n}} ,v_{\varepsilon_{n}}
  ,m_{\varepsilon_{n}} ,c_{\varepsilon_{n}} )$ solution to (\ref{eq:
  nondimensional travelling wave}) converges in $H^{1} ( \mathbbm{R} ) \times
  H^{1} ( \mathbbm{R} ) \times H^{1} ( \mathbbm{R} ) \times \mathbbm{R}$ to $(
  u,v,m,c )$, solution to the problem
  \begin{equation}
    \begin{array}{rcl}
      -u'' +c u' & = & c  \delta_{x= \bar{x}} ,\\
      - \Lambda v'' +c v' & = & -c  \delta_{x= \bar{x}} -c n_{0} m'  _{}
      \hspace{1em} \tmop{on}  \mathbbm{R},\\
      c m' & = & - \phi ( u,m ) \hspace{1em} \tmop{on}  \mathbbm{R},\\
      u ( - \infty ) =0, & u ( 0 ) = \theta , & u ( + \infty ) =1,\\
      v ( - \infty ) =v_{u} , &  & v ( + \infty ) =0,\\
      m ( - \infty ) =m_{u} . &  & 
    \end{array} \label{eq: limiting system in the HEA limit}
  \end{equation}
  Here $\bar{x} =- \log \theta_{i} /c$.
\end{theorem}

For any given $c>0$, the corresponding solution $( u,v,m )$ is uniquely
determined and one can find explicit expressions for all the unknowns,
provided the expression for $\phi$ can be explicitly integrated. The complete
characterisation of the system in the HAE limit is therefore dictated by the
exact value of the velocity $c>0$ of the travelling wave. The latter is
determined thanks to the analysis of the internal combustion layer, also
proved in Section \ref{SS: characterisation of the limiting profiles}:

\begin{theorem}
  \label{T: rigorous internal layer analysis}{\dueto{internal combustion layer
  analysis}}Let $( u,v,m,c )$ the limiting profile in the high activation
  energy limit. The velocity $c$ is given by
  \[ c =  \min \left\{ \sqrt{2 \mu / \Lambda} ,c_{\star} ( m_{u} ) \right\}, \]
  where $c_{\star} ( m_{u} )$ is the unique velocity such that the droplets
  finish vaporising exactly at the position of the reaction zone.
\end{theorem}

This result provides a rigorous justification of the existence of a
vaporisation controlled regime, i.e. a regime where $c=c_{\star} ( m_{u} ) <
\sqrt{2 \mu / \Lambda}$ (see {\cite{SuaNicHal2001}}). This regime appears when
the time for complete vaporisation of the droplets is large enough, precisely
when $m_{u} >m_{u}^{\star}$, where $m_{u}^{\star}$ is defined by $c_{\star} (
m_{u}^{\star} ) = \sqrt{2 \mu / \Lambda}$.

Section \ref{SS: internal layer in the vaporisation limited regime} provides
more details about the internal combustion layer in the vaporisation
controlled regime. An upper estimate of the size of the region where combustion and
vaporisation overlap is given in Theorem \ref{T: estimate of the overlapping
region}, followed by an application to the so--called $d^{2}$--law for
vaporisation.

Finally, Section \ref{S: extensions} presents some immediate applications of
our analysis to the case of polydisperse sprays, fast vaporisation, or radial
geometry. We refer to {\cite{Ber2003}} for a more detailed exposition of those
results.

\subsubsection*{Sketch of typical profiles.}The pictures below sketch
typical profiles for spray travelling waves moving from the right to the left.
The temperature profile $u$ is in red, the gaseous reactant profile $v$ is in
blue and the liquid phase density profile $n_{0} m$ is in green.
Recall that $u ( - \infty ) =0$, $v ( - \infty ) =v_{u}$, $m ( - \infty )
=m_{u}$ and $u ( + \infty ) =1=v_{u} +n_{0} m_{u}$. We set $n_{0} =1$.

\begin{figure}[ht]
\centering
\hfill
\begin{minipage}{6cm}
  \centering
  \includegraphics[width=6cm]{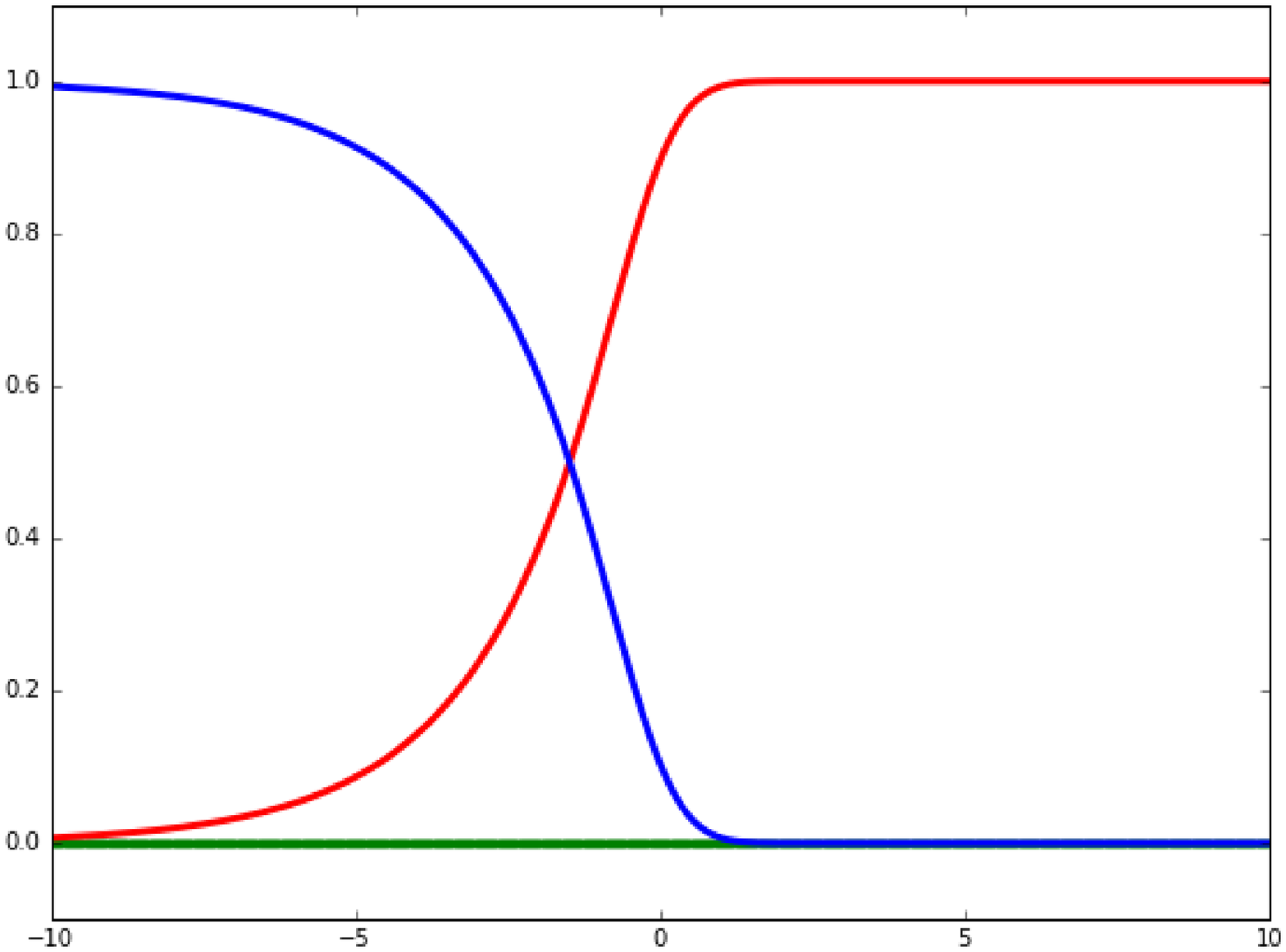}
  \captionof{figure}{}
  \label{fig: gasflame}
\end{minipage}
\hfill
\begin{minipage}{6cm}
  \centering
  \includegraphics[width=6cm]{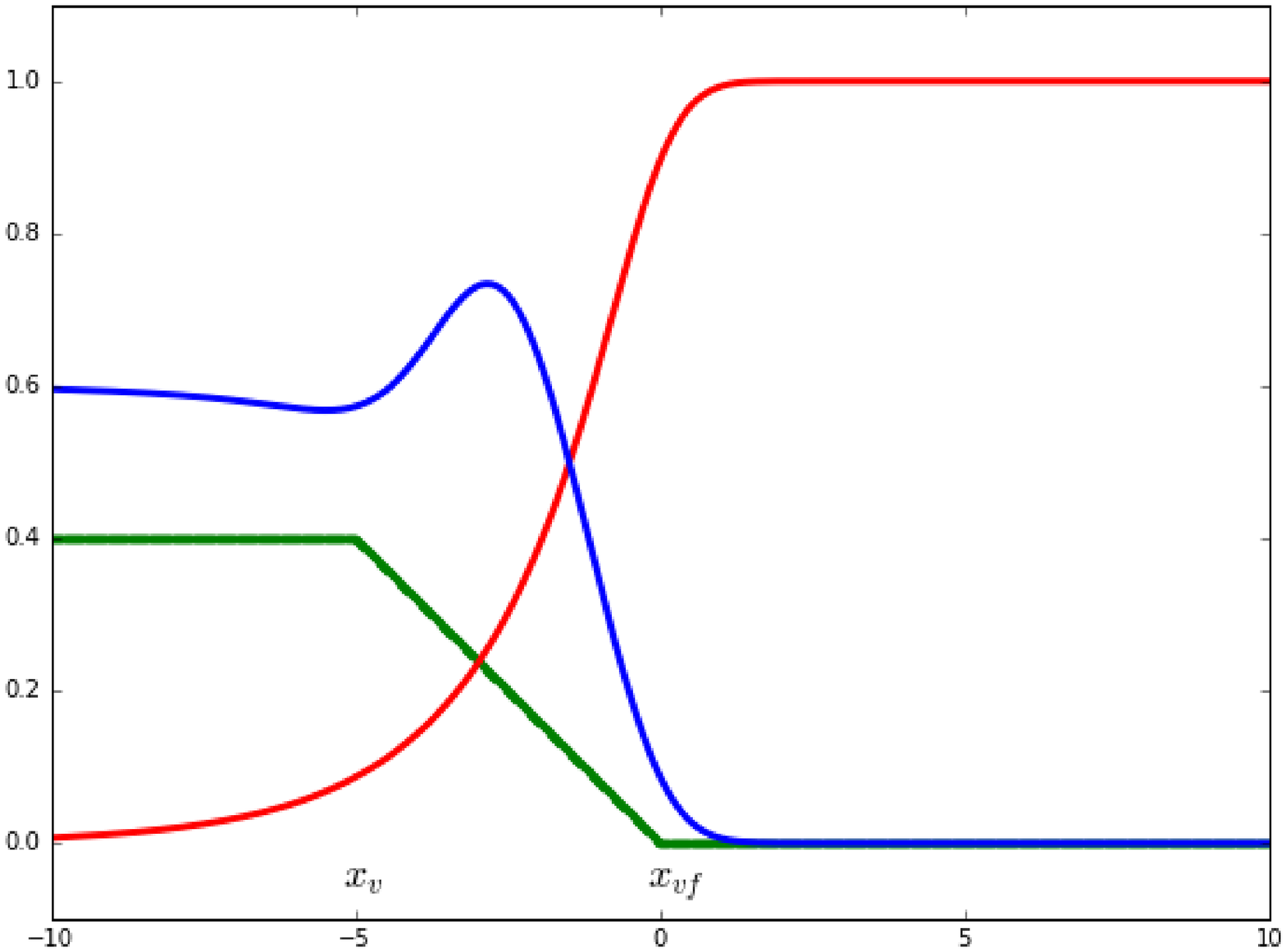}
  \captionof{figure}{}
  \label{fig: sprayflame}
\end{minipage}
\hfill\mbox{}
\end{figure}

Figure \ref{fig: gasflame} shows typical profiles for gas flames in the
absence of droplets, for a Lewis number equal to unity. In Figure \ref{fig:
sprayflame}, droplets are present in the fresh gas with mass density $n_{0}
m_{u} =0.4$. We therefore set $v_{u} =0.6$ so that the temperature in the
burnt gas remains equal to unity. Notice that the gaseous reactant profile is
no longer monotoneous, due to the vaporisation of the liquid phase.

\begin{figure}[ht]
\centering
\hfill
\begin{minipage}{6cm}
  \centering
  \includegraphics[width=6cm]{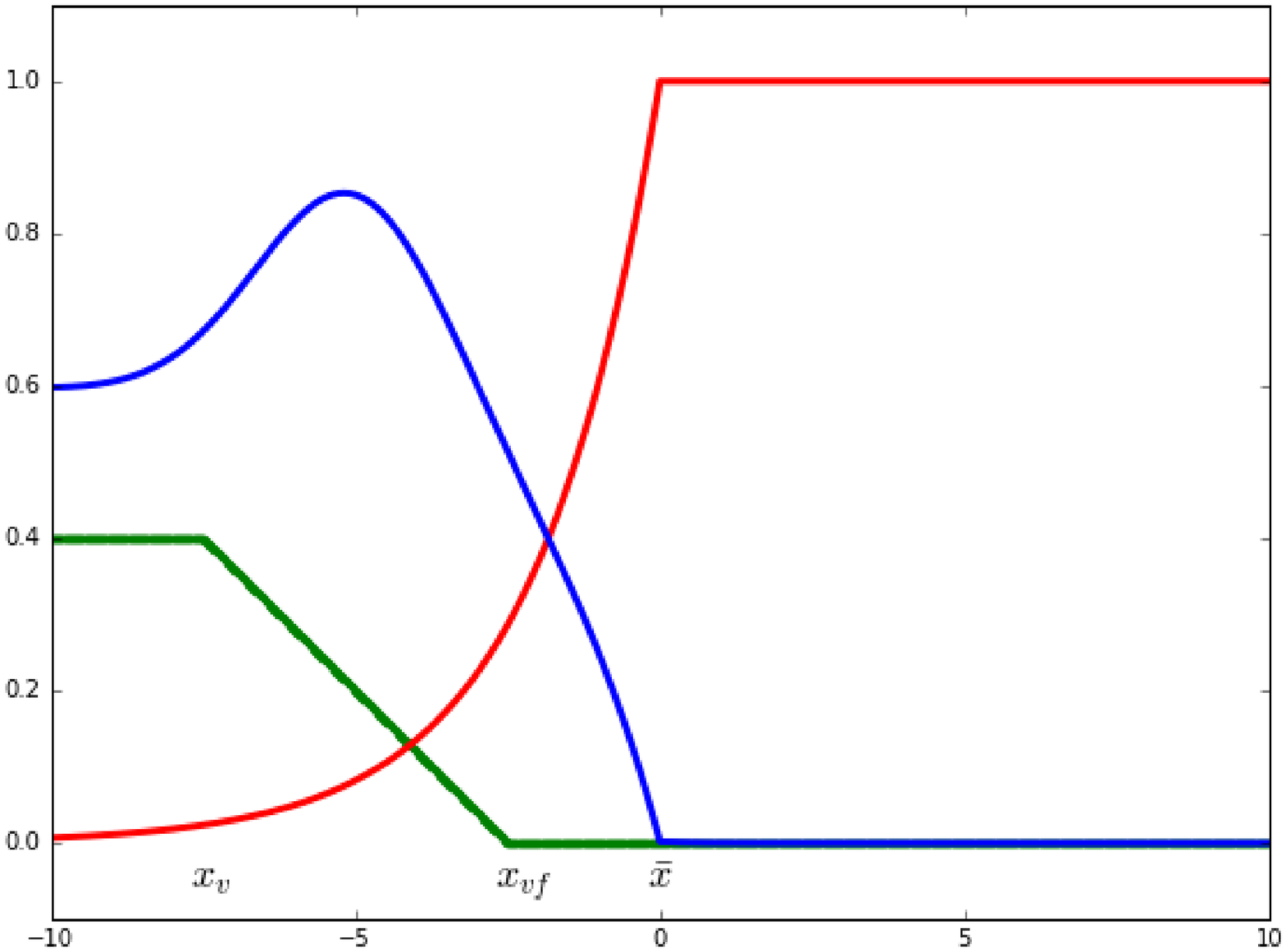}
  \captionof{figure}{}
  \label{fig: HAE}
\end{minipage}
\hfill
\begin{minipage}{6cm}
  \centering
  \includegraphics[width=6cm]{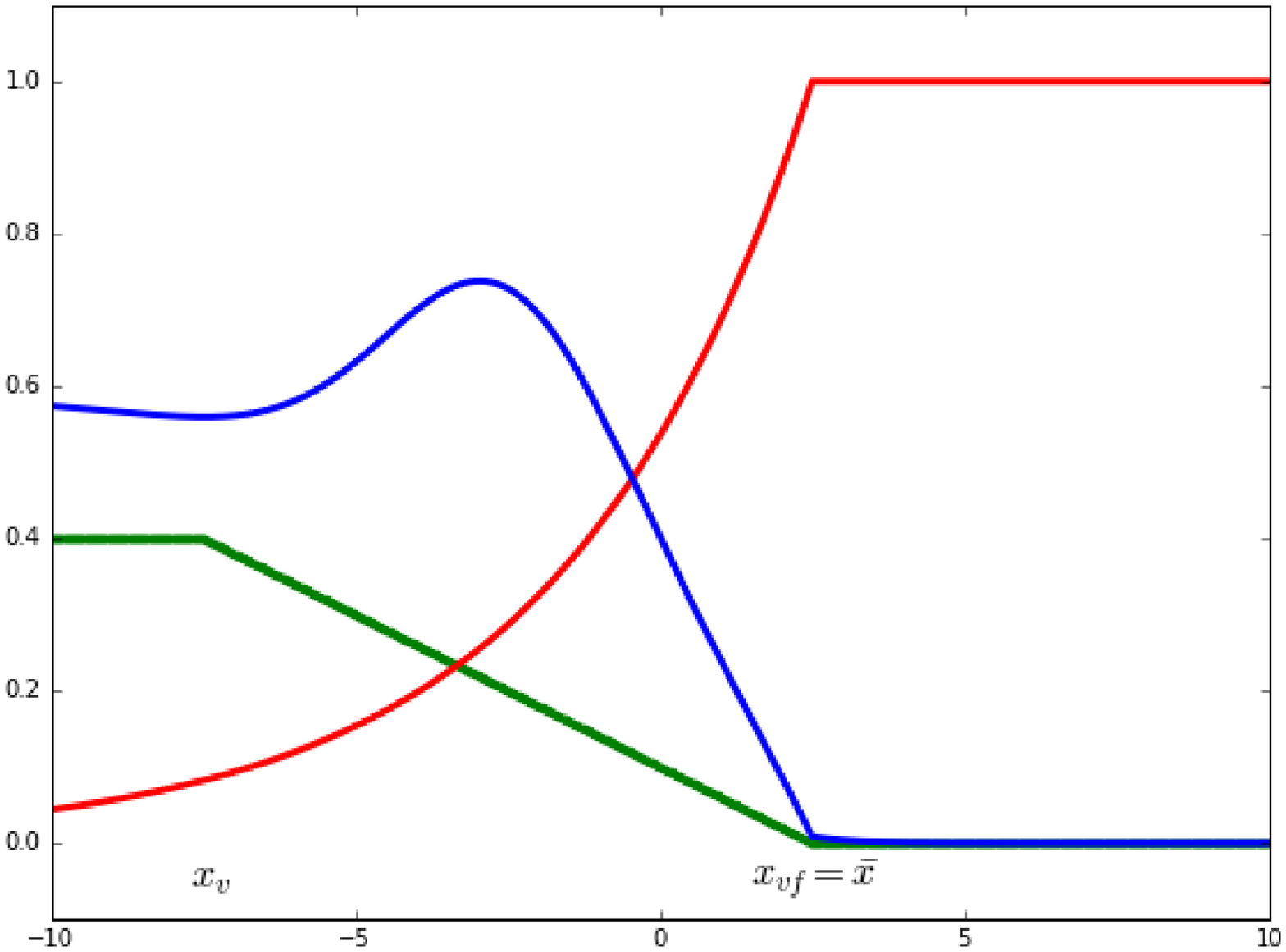}
  \captionof{figure}{}
  \label{fig: HAElimited}
\end{minipage}
\hfill\mbox{}
\end{figure}

Figures \ref{fig: HAE} \& \ref{fig: HAElimited} are typical profiles in the
high activation energy limit. The reaction zone is reduced to a point located
at $\bar{x} =0$ in Figure \ref{fig: HAE} and at $\bar{x} =2.5$ in Figure~\ref{fig:
HAElimited}. In both cases vaporisation starts at $x_{v} =-7.5$. However, in
the first case the vaporisation front, i.e. the point $x_{v f}$ where
vaporisation ends, is located {\tmstrong{before}} the reaction front, at
$x_{v} =-2.5<0$. Our analysis will show that in such a situation, the velocity
of the spray flame in the HAE limit is equal to that of the purely gaseous
flames with same temperature in the burnt gas. In the second case, the
vaporisation ends at $x_{v f} =2.5$. Shall the temperature profile remain the
same as in Figure \ref{fig: HAE}, the vaporisation process would still occur
{\tmstrong{after}} the reaction zone. Our analysis will show that this is an
\tmstrong{impossibility} and that the preheating zone ``has to'' stretch so that the
reaction front $\bar{x}=0$ from Figure {\ref{fig: HAE}} is ``pushed'' to the right
until it coincides with the vaporisation front $\bar{x} =x_{ v f}=2.5$ as in Figure
\ref{fig: HAElimited}. The stretching of the preheating zone is synonimous
with a decrease of the velocity of the flame. This is the so--called vaporisation controlled regime
(see {\cite{SuaNicHal2001}}).

\section{Existence of travelling waves}\label{S: Existence of travelling
waves}

Let $[ -a,a ]$ an interval of $\mathbbm{R}$ for some $a>0$. Let $X_{a} =
\mathcal{C}^{1} ( [ -a,a ] ) \times \mathcal{C}^{1} ( [ -a,a ] ) \times
\mathcal{C}^{0} ( [ -a,a ] ) \times \mathbbm{R}$. We consider system (\ref{eq:
nondimensional travelling wave}) restricted to the interval $[ -a,a ]$ with
Dirichlet boundary conditions at $x=a$ and flux conditions at $x=-a$,
\begin{equation}
  \begin{array}{rcl}
    -u'' +c u' & = & f ( u ) v \hspace{1em} \tmop{on}   ( -a,a ) ,\\
    - \Lambda v'' +c v' & = & -f ( u ) v-c n_{0} m'  _{} \hspace{1em}
    \tmop{on}   ( -a,a ) ,\\
    c m' & = & - \phi ( u,m ) \hspace{1em} \tmop{on}   ( -a,a ) ,\\
    -u' ( -a ) +c u ( -a ) =0, &  & u ( a ) =1,\\
    - \Lambda v' ( -a ) +c v ( -a ) =c v_{u} , &  & v ( a ) =0,\\
    m ( -a ) =m_{u} , &  & \\
    u ( 0 ) = \theta_{i} . &  & 
  \end{array} \label{eq: system on a bounded domain}
\end{equation}
The last equation allows one to break the translation invariance in the limit
where $a$ goes to infinity. We first prove \tmtextit{a priori} estimates for
solutions $( u,v,m,c ) \in X_{a}$. This allows us to establish the
existence of solutions in a bounded domain by using topological degree
arguments. Finally, the existence is extended to the real line.

\subsection{Preliminary estimates}\label{SS: preliminary estimates}The main
result of this section consists of \tmtextit{a priori} estimates for solutions
to problem (\ref{eq: system on a bounded domain}) above that are uniform with
respect to the size $a$ of the domain. We prove successively the following
three results:

\begin{proposition}
  \label{P: qualitative properties}{\dueto{qualitative properties}}Let $(
  u,v,m,c ) \in X_{a}$ a solution to (\ref{eq: system on a bounded domain})
  with $c \geqslant 0$. We have,
  \begin{eqnarray}
    & c>0, & \\
    & 0<u \leqslant 1, \hspace{1em} 0 \leqslant v<1, \hspace{1em} 0 \leqslant
    m \leqslant m_{u} \hspace{1em} \tmop{on}   [ -a,a ], & \\
    & 0<u' \leqslant c, \hspace{1em} -c \leqslant \Lambda v' \leqslant c (
    1+n_{0} m_{u} ) , \hspace{1em} m' \leqslant 0, \hspace{1em} \tmop{on}   [
    -a,a ]. & 
  \end{eqnarray}
\end{proposition}

\begin{proposition}
  \label{P: a priori estimates}{\dueto{\tmtextit{a priori} estimates}}Let $(
  u,v,m,c )$ a solution to (\ref{eq: system on a bounded domain}) with $c
  \geqslant 0$. We have
  \begin{eqnarray*}
    \beta_{1} ( \Lambda^{-1} ) ( 1-u ( x ) ) -n_{0} m ( x ) & \leqslant v ( x
    ) \leqslant & \beta_{2} ( \Lambda ) ( 1-u ( x ) )  ,
  \end{eqnarray*}
  where $\beta_{1} ( \Lambda ) \assign \min ( 1, \Lambda^{-1} )$ and
  $\beta_{2} ( \Lambda ) \assign \max ( 1, \Lambda^{-1} )$.
\end{proposition}

In the following we note, for all $\theta_{i} \leqslant s \leqslant 1$,
\[ G ( s ) = \int_{\theta_{i}}^{s} f ( s ) ( 1-s ) \mathd s. \]
\begin{proposition}
  \label{P: bounds on the velocity}{\dueto{bounds on the velocity}}Let
  $\theta^{\star} \assign ( \theta_{v} + \theta_{i} ) /2$, $c_{\star} =
  \sqrt{\tau ( m_{u} , \theta^{\star} ) / \log ( \theta_{i} / \theta^{\star}
  )}$ and $a_{\star} = \log ( \theta_{i} / \theta_{v} ) /c_{\star}$. For all
  $a \geqslant a_{\star}$, the solution $( u,v,m,c )$ to problem (\ref{eq:
  system on a bounded domain}) with $c \geqslant 0$ obeys
  \begin{equation}
    \min ( c_{\star} ,c_{1} ) \leqslant c \leqslant {\color{red} } ( c_{2} ( a
    ) ,c_{3} ) , \label{eq: bounds on the velocity}
  \end{equation}
  where we introduced  $c_{1} = \sqrt{2 \beta_{1} ( \Lambda ) G ( 1 )} / \theta_{i}$, $c_{2} (
  a ) = \sqrt{2 \beta_{2} ( \Lambda ) ( G ( 1 ) +u' ( a )^{2} )} /
  \theta_{i}$, $M_{f} = \sup_{s \in [ 0,1 ]} f ( s )$, 
  $c_{3} = \max \left( M_{f} , \log ( 4/ \theta_{i} ) /a_{\star}
  ,2 \sqrt{M_{f}/ \theta_{i}} \right)$.
\end{proposition}

\begin{proof*}{Proof of Proposition \ref{P: qualitative properties}
{\dueto{qualitative properties}}}
  The presence of vaporising fuel in the model does not allow one to recover
  all the qualitative properties found in the gaseous case. For example, an important
  difference is the non--monotonicity of the reactant profile $v ( x )$.
  Also, we introduce below auxiliary functions (such as $v_{\ast}$, $\bar{w}$ or $\bar{v}$)
  that are specific to the analysis of spray flames,
  and specific steps such as the proof that $v\geqslant 0$  or the localisation
  of the vaporisation region.  
  
  \subsubsection*{Proof that $c>0$.}Recall that we assumed $c \geqslant 0$. Assume by
  contradiction that $c=0$. Then the vaporisation term $c n_{0} m'$ disappears
  in the equation of the reactant $v$. The function $y \assign u+ \Lambda v-1$
  therefore solves $y'' =0$ on $( -a,a )$, $y ( -a ) =y ( a ) =0$, hence $y=0$
  on $[ -a,a ]$. It follows that $( u,c )$ solves the nonlinear
  Dirichlet--Neumann boundary value problem
  \[ \begin{array}{l}
       -u'' - \displaystyle\frac{1}{\Lambda} f ( u ) ( 1-u ) =0 \hspace{1em} \tmop{on}   (
       -a,a ) ,\medskip\\
       u' ( -a ) =0, \hspace{1em} u ( 0 ) = \theta , \hspace{1em} u ( a ) =1.
     \end{array} \]
  We then observe that $u \leqslant 1$ on $[ -a,a ]$ otherwise the maximum $u
  ( x_{0} ) >1> \theta$ of $u$ would be attained at $x_{0} \in [ -a,a )$ with
  $f ( u ( x_{0} ) ) >0$. The equation above would imply $u'' ( x_{0} ) >0$ in
  contradiction with $x_{0}$ being a maximum. It follows that $w \assign u'$
  obeys $w' \leqslant 0$ and $w ( -a ) =0$, hence $w \leqslant 0$ and $u$
  decreasing on $[ -a,a ]$; a contradiction with $u ( 0 ) = \theta <1=u ( a
  )$.
  
  \subsubsection*{Proof that $v \geqslant 0$.}Given a solution $( u,v,m,c ) \in X_{a}$,
  introduce $v_{\ast}$ the auxiliary reactant profile solution to the linear
  elliptic boundary value problem
  \[ \begin{array}{l}
      \mathcal{L} v_{\ast} \assign - \Lambda v_{\ast}'' +c v_{\ast}' -f ( u )
      v_{\ast} =0 \hspace{1em} \tmop{on}   ( -a,a ) , \medskip\\
      \mathcal{N} v_{\ast} ( -a ) \assign - \Lambda v_{\ast}' ( -a ) +c v_{\ast} ( -a ) =0,
      \hspace{1em} v_{\ast} ( a ) =0. 
      \end{array}\]
  The positive function $\phi \equiv 1$ obeys $\mathcal{L} \phi \geqslant 0$
  on $( -a,a )$, $\mathcal{N} \phi ( -a ) =c>0$. Therefore the problem has a
  unique solution, obviously $v_{\ast} \equiv 0$, and the generalized maximum
  principle implies, since $\mathcal{L} v=-c n_{0} m' \geqslant 0$,
  $\mathcal{N} v ( -a ) =c v_{u} \geqslant 0$ and $v ( a ) \geqslant 0$, that
  $v \geqslant v_{\ast} \geqslant 0$ (see e.g. {\cite[Chapter 2]{ProWei1984a}}).
  
  \subsubsection*{Bounds for $u$, $v$ and $m$.}The bound $0 \leqslant m ( x )
  \leqslant m_{u}$ is trivial. Let $w \assign - \Lambda v' +c v+c n_{0} m$. We
  have $w' =-f ( u ) v \leqslant 0$, therefore $c=w ( -a ) \geqslant w ( x )
  \geqslant w ( a ) =- \Lambda v' ( a ) +c n_{0} m ( a ) \geqslant 0$. We used
  $v' ( a ) \leqslant 0$, a consequence of $v \geqslant 0$ and $v ( a ) =0$.
  Integrating the equation $-u'' +c u' =-w'$ on $[ -a,x ]$ yields $-u' ( x )
  +c u ( x ) =-e^{c x} ( e^{-c x} u )' =c-w ( x )$. Therefore $-c e^{-c x}
  \leqslant ( e^{-c x} u )' \leqslant 0$ and integrating now on $[ x,a ]$
  yields exactly $0 \leqslant u ( x ) \leqslant 1$.
  
  In order to bound $v$ from above, we introduce the auxiliary function
  $\bar{v}$, solving
  \begin{equation}
     \begin{array}{l}
     - \Lambda \bar{v}'' +c \bar{v}' =-f ( u ) v \hspace{1em} \tmop{on}
    ( -a,a ),\medskip\\
     - \Lambda \bar{v}' ( -a ) +c \bar{v}
    ( -a ) =c, \hspace{1em} \bar{v} ( a ) =0. \label{eq: vbar}
    \end{array}
  \end{equation}
  Then let $\bar{w} \assign \bar{v} -v$, solution to $- \Lambda \bar{w}'' +c
  \bar{w}' =c n_{0} m'$ on $( -a,a )$, $- \Lambda \bar{w} ( -a ) +c \bar{w} (
  -a ) =c n_{0} m_{u}$ and $w ( a ) =0$. Integrating over $[ -a,x ]$ yields \
  $- \Lambda e^{c x/ \Lambda} ( e^{-c x/ \Lambda} \bar{w} )' =c n_{0} m ( x
  )$, and integrating now on $[ x,a ]$ yields
  \begin{equation}
    \Lambda \bar{w} ( x ) =c n_{0} \int_{x}^{a} e^{-c ( s-x ) / \Lambda}  m (
    s ) \mathd s. \label{eq: wbar}
  \end{equation}
  obviously $\bar{w} \geqslant 0$, that is $v \leqslant \bar{v}$ and it
  remains to bound $\bar{v}$ from above. For that we observe that $- \Lambda
  \bar{v}'' +c \bar{v}' =+w'$. Integrating on $[ -a,x ]$ yields $( -c/ \Lambda
  ) e^{-c x/ \Lambda} \leqslant ( e^{-c x/ \Lambda} \bar{v} )' \leqslant 0$,
  and integrating now on $[ x,a ]$ yields exactly \ $0 \leqslant \bar{v}
  \leqslant 1$.
  
  \subsubsection*{Bounds for $u'$ and $v'$.}We have from the previous step $u' ( x )
  =c ( u ( x ) -1 ) +w ( x ) \leqslant w ( x ) \leqslant c$ since $u \leqslant
  1$ and $w \leqslant c$. Also the equation for $u$ rewrites $-e^{c x} ( e^{-c
  x} u' )' =f ( u ) v \geqslant 0$ and integrating from $x$ to $a$ yields
  $e^{-c a} u' ( a ) \leqslant e^{-c x} u' ( x )$. Integrating the equation $-
  ( u+ \Lambda \bar{v} )'' +c ( u+ \bar{v} )' =0$ on $( -a,a )$ allows us to
  estimate $u' ( a ) =- \Lambda \bar{v}' ( a ) >0$, since $\bar{v} \geqslant v
  \geqslant 0$, $\bar{v} ( a ) =0$ and $\bar{v}$ is not identically zero. The
  monotonicity $u' >0$ follows. Let again $w \assign - \Lambda v' +c v+c n_{0}
  m$. The bounds on $0 \leqslant w \leqslant c$, $0 \leqslant v \leqslant 1$
  and $0 \leqslant m \leqslant m_{u}$ easily imply $-c \leqslant v' \leqslant
  c ( 1+n_{0} m_{u} )$.
  
  This concludes the proof of of Proposition \ref{P: qualitative properties}
  {\dueto{qualitative properties}}.
\end{proof*}

\begin{proof*}{Proof of Proposition \ref{P: a priori estimates}
{\dueto{\tmtextit{a priori} estimates}}}
  In order to obtain a pointwise comparison of $v$ with $( 1-u )$, we take a
  detour and compare rather $\bar{v}$ with $( 1-u )$, where $0 \leqslant
  \bar{v} \leqslant 1$ is the auxiliairy function defined in (\ref{eq: vbar}).
  Observe first that $\bar{v}$ is nonincreasing. Indeed from (\ref{eq: vbar})
  we have $( e^{-c x/ \Lambda} \bar{v}' )' \geqslant 0$, and integrating on $[
  x,a ]$ yields $e^{-c x/ \Lambda} \bar{v}' ( x ) \leqslant e^{-c a/ \Lambda}
  \bar{v}' ( a ) \leqslant 0$.
  
  Following {\cite{BerNicSch1985}}, introduce now the two auxiliary functions
  $z \assign u+ \bar{v} -1$ and $y \assign u+ \Lambda \bar{v} -1$ and the four
  relations
  \begin{eqnarray*}
    -z'' +c z' & = & ( \Lambda -1 ) \bar{v}'' ,\\
    - \Lambda z'' +c z' & = & ( 1- \Lambda ) u'' ,\\
    -y'' +c y' & = & ( \Lambda -1 ) c \bar{v}' ,\\
    - \Lambda y'' +c y' & = & ( 1- \Lambda ) \bar{c} u' .
  \end{eqnarray*}
  Integrating these equations successively on $[ -a,x ]$ and $[ x,a ]$ yields
  \begin{eqnarray*}
    z ( x ) & = & ( \Lambda -1 ) \int_{x}^{a} e^{-c ( s-x )}   \bar{v}' ( s ) 
    \mathd s,\\
    z ( x ) & = & \frac{1- \Lambda}{\Lambda} \int_{x}^{a} e^{-c ( s-x ) /
    \Lambda}  u' ( s )   \mathd s,\\
    y ( x ) & = & ( 1- \Lambda ) c \int_{x}^{a} e^{-c ( s-x )} \bar{v} ( s )  
    \mathd s,\\
    y ( x ) & = & \frac{\Lambda -1}{\Lambda} c \int_{x}^{a} e^{-c ( s-x ) /
    \Lambda} ( 1-u ) ( s )   \mathd s.
  \end{eqnarray*}
  Since $\bar{v}$ and $( 1-u )$ are nonincreasing, it is not difficult to
  obtain
  \[ \begin{array}{l}
       | z ( x ) | \leqslant | \Lambda -1 |   \bar{v} ( x ) , \hspace{2em} | z
       ( x ) | \leqslant \left\vert \displaystyle\frac{\Lambda -1}{\Lambda} \right\vert   ( 1-u ( x
       ) ),\medskip\\
       | y ( x ) | \leqslant | \Lambda -1 |   \bar{v} ( x ) , \hspace{2em} | y
       ( x ) | \leqslant | \Lambda -1 |   ( 1-u ( x ) ) ,
     \end{array} \]
  which implies
  \[ \begin{array}{ll}
       \displaystyle\frac{1}{\Lambda} ( 1-u ( x ) ) \leqslant \bar{v} ( x ) \leqslant ( 1-u
       ( x ) ) & \hspace{2em} \tmop{if} \hspace{1em} \Lambda >1\medskip\\
       ( 1-u ( x ) ) \leqslant \bar{v} ( x ) \leqslant \displaystyle\frac{1}{\Lambda} ( 1-u
       ( x ) ) & \hspace{2em} \tmop{if} \hspace{1em} 0< \Lambda <1.
     \end{array} \]
  It remains to compare precisely $\bar{v}$ and $v= \bar{v} - \bar{w}$. But
  since $m$ is nonincreasing nonnegative, the expression (\ref{eq: wbar})
  ensures pointwise $0 \leqslant \bar{w} ( x ) \leqslant n_{0} m ( x )$. The
  result announced follows.

  This concludes the proof of Proposition \ref{P: a priori estimates}
  {\dueto{\tmtextit{a priori} estimates}}.
\end{proof*}

\begin{proof*}{Proof of Proposition \ref{P: bounds on the velocity}
{\dueto{bounds on the velocity}}}
  The nonlinear eigenvalue $c$ is estimated thanks to energy estimates.
  
  \subsubsection*{Estimates on $[ 0,a ]$.}Integrating successively the equation $f (
  u ) v=-u'' +c u'$ against $1$, $u$ and $u'$ on $[ 0,a ]$ yields
  \begin{eqnarray*}
    \int_{0}^{a} f ( u ) v \cdot 1 & = & c-u' ( a ), \phantom{+ \frac{c}{2} ( 1+
    \theta_{i} )^{2} + \int_{0}^{a} ( u' )^{2}} \mbox{} \hspace{2.2cm} ( \star )\\
    \int_{0}^{a} f ( u ) v \cdot u & = & \phantom{c}-u' ( a ) + \frac{c}{2} ( 1+
    \theta_{i} )^{2} + \int_{0}^{a} ( u' )^{2}, \hspace{2cm} ( \star \star )\\
    \int_{0}^{a} f ( u ) v \cdot u' & = & \phantom{c} - \frac{1}{2} u' ( a )^{2} +
    \frac{1}{2}  c^{2} \theta_{i}^{2} +c \int_{0}^{a} ( u' )^{2}. \hspace{2cm} (
    \star \star \star )
  \end{eqnarray*}
  where we used $u ( 0 ) = \theta_{i}$, $u' ( 0 ) =c \theta_{i}$. The
  combination $( \star \star \star ) +c [ ( \star ) - ( \star \star ) ]$ reads
  \begin{eqnarray*}
    ( \star \star \star \star ) & \assign & \int_{0}^{a} f ( u ) v \cdot u' +c
    \int_{0}^{a} f ( u ) v \cdot ( 1-u )\\
    & = & - \frac{1}{2} u' ( a )^{2} + \frac{c^{2}}{2}   \theta_{i}^{2} .
  \end{eqnarray*}
  Since $( 1-u ) \geqslant 0$, the two equations $( \star \star \star )$ and
  $( \star \star \star \star )$ imply
  \[ \int_{0}^{a} f ( u ) v \cdot u' \leqslant \frac{1}{2}  c^{2}
     \theta_{i}^{2} \leqslant \int_{0}^{a} f ( u ) v \cdot u' + \frac{1}{2} u'
     ( a )^{2}. \]
  It remains to estimate the contribution of the reaction term in the
  inequalities above.
  
  \subsubsection*{First upper bound for $c$.}From Proposition \ref{P: a priori
  estimates}, we have $v ( x ) \leqslant \beta_{2} ( \Lambda ) ( 1-u ( x ) )$
  on $[ 0,a ]$, hence
  \[ \int_{0}^{a} f ( u ) v \cdot u' \leqslant \beta_{2} ( \Lambda )
     \int_{0}^{a} f ( u ( x ) ) ( 1-u ( x ) ) \cdot u' ( x ) \mathd x=
     \beta_{2} ( \Lambda ) G ( 1 ), \]
  which provides the value of $c_{2} ( a )$ in Proposition \ref{P: bounds on
  the velocity} such that
  \[ \frac{1}{2}  c_{2} ( a )^{2}   \theta_{i} \assign \beta_{2} ( \Lambda ) G
     ( 1 ) + \frac{1}{2} u' ( a )^{2} . \]
  \subsubsection*{Localisation of the vaporisation region.}We claim that for small
  enough velocities, the vaporisation region does not intersect the combustion
  zone. Recall that $\theta_{v} < \theta_{i}$ so that $\theta^{\star} \assign
  ( \theta_{v} + \theta_{i} ) /2< \theta_{i}$. Let $\tau ( m_{u} ,
  \theta^{\star} )$ the interval of time needed for complete vaporisation of a
  droplet of size $m_{u}$ at constant temperature $\theta^{\star}$. From the
  hypothesis on the vaporisation term, we have $\tau ( m_{u} , \theta^{\star}
  ) < \infty$. Given $c>0$, let $x_{v}$ such that $u ( x_{v} ) = \theta_{v}$,
  $x^{\star}$ such that $u ( x^{\star} ) = \theta^{\star}$ and $x_{i} =0$ such
  that $u ( x_{i} ) = \theta_{i}$. Since $u ( x ) = \theta_{i} e^{c x}$ on $[
  -a \nocomma ,0 ]$, we have for any given $c$ and $a$ large enough, $-a
  \leqslant x_{v} <x^{\star} <x_{i}$ and $x_{i} -x^{\star} = \log ( \theta_{i}
  / \theta^{\star} ) /c$. The time $\tau^{\star}$ spent by the droplets
  advected at velocity $c$ inside the interval $[ x^{\star} ,x_{i} ]$ is
  $\tau^{\star} = ( x_{i} -x^{\star} ) /c= \log ( \theta_{i} / \theta^{\star}
  ) /c^{2}$. However, since $m ( x^{\star} ) \leqslant m_{u}$, and $u
  \geqslant \theta^{\star}$ on $[ x^{\star} ,x_{i} ] \nocomma$, the
  monotonicity properties of the vaporisation law imply that if $\tau^{\star}
  \geqslant \tau ( m_{u} , \theta^{\star} )$ hence complete vaporisation
  occurs inside the intervall $[ x^{\star} ,x_{i} ]$. To prove our claim, it
  suffices to take $a \geqslant a_{\star}$ and $c \leqslant c_{\star}$, with
  \[ c_{\star} = \sqrt{\tau ( m_{u} , \theta^{\star} ) / \log ( \theta_{i} /
     \theta^{\star} )} , \hspace{1em} a_{\star} = \log ( \theta_{i} /
     \theta_{v} ) /c_{\star} . \]

  \subsubsection*{Lower bound for $c$.}Take $a \geqslant a_{\star}$. If $c \geqslant
  c_{\star}$ there is nothing to prove. Otherwise, we have $c \leqslant
  c_{\star}$, $m ( x ) =0$ on $[ 0,a ]$ so that from Proposition \ref{P: a
  priori estimates} we have $v ( x ) \geqslant \beta_{1} ( \Lambda ) ( 1-u ( x
  ) )$, therefore
  \[ \int_{0}^{a} f ( u ) v \cdot u' \geqslant \beta_{1} ( \Lambda )
     \int_{0}^{a} f ( u ( x ) ) ( 1-u ( x ) ) \cdot u' ( x ) \mathd x=
     \beta_{1} ( \Lambda ) G ( 1 ),
  \]
  which allows us to define the value $c_1$ in Proposition \ref{P: bounds on the
  velocity} thanks to the relation $c_{1}^{2}   \theta_{i}/2 \assign \beta_{1} ( \Lambda ) G ( 1 )$.
  
  \subsubsection*{Second upper bound for $c$.}Let $a_{\star}$ defined above. Let $M
  \assign \sup_{s \in [ 0,1 ]} f ( s )$. Let $\bar{u}$ be the solution to $-
  \bar{u}'' +c \bar{u}' =M  \mathbbm{1}_{] 0,a [}$ on $] -a,a [$, $- \bar{u}'
  ( -a ) +c \bar{u} ( -a ) =0$ and $\bar{u} ( a ) =1$. Since $0 \leqslant v
  \leqslant 1$, the maximum principle asserts that $u \leqslant \bar{u}$ on $[
  -a,a ]$, in particular $\theta_{i} \leqslant \bar{\theta} \assign \bar{u} (
  0 )$. Solving explicitly for $\bar{u}$ yields
  \[ \bar{u} ( 0 ) =e^{-c a} ( 1-M/c ) + ( 1-e^{-c a} ) M/c^{2} . \]
  Assume $a \geqslant a_{\star} >0$, then the right hand side tends towards
  zero as $c$ goes to infinity. More precisely, if both $c \geqslant M$,
  $e^{-c a_{\star}} \leqslant \theta_{i} /4$ and $M/c^{2} \leqslant \theta_{i}
  /4$ one has $\bar{u} ( 0 ) = \bar{\theta} \leqslant \theta_{i} /2$, a
  contradiction. It follows that
  \[ c_{3} \assign \max \left( M, \log ( 4/ \theta_{i} ) /a_{\star} ,2
     \sqrt{M/ \theta_{i}} \right) \]
  is an upper--bound for $c$.

  This concludes the proof of Proposition \ref{P: bounds on the velocity}
  {\dueto{bounds on the velocity}}.
\end{proof*}

\subsection{Existence of a solution on bounded domains}\label{SS: existence of
a solution on bounded domains}

The proof of the existence of a topological degree relies on the existence of
the continuous map $K_{\tau}$ below that is not continuous if the reaction or
vaporisation terms have discontinuities at $\theta_{i}$ or $\theta_{v}$
respectively. Following
{\tmname{Berestycki}}--{\tmname{Nicolaenko}}--{\tmname{Scheurer}} (see
{\cite[footnote p. 1225]{BerNicSch1985}}), we assume below that $f$ and $g$
are continuous. The case where a discontinuity is present can be attained by a
standard smoothing procedure.

The main point when dealing with a spray flame model is to correctly account
for the vaporisation terms in the homotopy argument in order to preserve
the  \tmtextit{a priori} estimates of the previous section.

\begin{proposition}
  \label{P: existence of a solution on bounded domains}{\dueto{existence of a
  solution on bounded domains}}There exists $a_{0} >0$ such that for $a
  \geqslant a_{0}$ system (\ref{eq: system on a bounded domain}) admits a
  solution in $X_{a} = \mathcal{C}^{1} ( [ -a,a ] ) \times \mathcal{C}^{1} ( [
  -a,a ] ) \times \mathcal{C}^{0} ( [ -a,a ] ) \times \mathbbm{R}$.
\end{proposition}

\begin{proof*}{Proof.}
  We use the {\tmname{Leray}}--{\tmname{Schauder}} topological degree
  argument. Let $\tau \in [ 0,1 ]$ the homotopy parameter and consider the
  solutions $( u_{\tau} ,v_{\tau} ,m_{\tau} ,c_{\tau} )$ to the new system
  \begin{equation}
    \begin{array}{rcl}
      -u'' +c u' & = & \tau \{ f ( u ) v \} \hspace{.5em} \tmop{on}   ( -a,a )
      ,\\
      - \Lambda v'' +c v' & = & \tau \{ -f ( u ) v+n_{0} \phi ( u,m ) \}  _{}
      \hspace{.5em} \tmop{on}   ( -a,a ) ,\\
      c m' & = & \tau \{ - \phi ( u,m ) \} \hspace{.5em} \tmop{on}   ( -a,a )
      ,\\
      -u' ( -a ) +c u ( -a ) &=&0, \hspace{1cm} u ( a ) =1,\\
      - \Lambda v' ( -a ) +c v ( -a ) &=&c v_{u} , \hspace{1cm} v ( a ) =0,\\
      m ( -a ) &=&m_{u} ,  \\
      c & = & u ( 0 ) - \theta_{i} + \tau c.
    \end{array} \label{eq: system with tau}
  \end{equation}
  This is system (\ref{eq: system on a bounded domain}) where the condition $u
  ( 0 ) = \theta_{i}$ has been rewritten as a fixed point for the velocity $c$
  and where the homotopy parmeter appears in the right--hand--side of the
  system. It is important to notice that this system can be obtained by
  replacing the vaporisation terms $f ( u )$ and $\phi ( u,m )$ in (\ref{eq:
  system on a bounded domain}) by $f_{\tau} \assign \tau f$ and $\phi_{\tau}
  \assign \tau \phi$. It follows that all \tmtextit{a priori} estimates from
  the previous sections hold.

  \subsubsection*{Fixed point formulation.}We look for solutions $( u_{\tau}
  ,v_{\tau} ,m_{\tau} ,c_{\tau} )$ defined as a fixed point of the mapping
  $K_{\tau} :X_{a} \rightarrow X_{a}$ that maps $( u,v,m,c )$ to $( \hat{u} ,
  \hat{v} , \hat{m} , \hat{c} )$ solution to the linear boundary value problem
  \begin{equation}
    \begin{array}{rcl}
      - \hat{u}'' +c  \hat{u}' & = & \tau \{ f ( u ) v \} \hspace{.5em}
      \tmop{on}   ( -a,a ) ,\\
      - \Lambda \hat{v}'' +c  \hat{v}' & = & \tau \{ -f ( u ) v+n_{0} \phi (
      u,m ) \}  _{} \hspace{.5em} \tmop{on}   ( -a,a ) ,\\
      c  \hat{m}' & = & \tau \{ - \phi ( u,m ) \} \hspace{.5em} \tmop{on}   (
      -a,a ) ,\\
      - \hat{u}' ( -a ) +c  \hat{u} ( -a ) &=&0, \hspace{1cm} \hat{u} ( a ) =1,\\
      - \Lambda \hat{v}' ( -a ) +c  \hat{v} ( -a ) &=&c v_{u} , \hspace{1cm} \hat{v} ( a
      ) =0,\\
      \hat{m} ( -a ) &=& m_{u} ,  \\
      \hat{c} & = & u ( 0 ) - \theta_{i} + \tau c.
    \end{array} \label{eq: system for the fixed point problem}
  \end{equation}
  Since $H^{2} ( ] -a,a [ )$ embeds compactly in $\mathcal{C}^{1} ( [ -a,a ]
  )$ and $H^{1} ( ] -a,a [ )$ embeds compactly in $\mathcal{C} ( [ -a,a ] )$,
  it follows that $K_{\tau}$ is a compact mapping and uniformly continuous
  with respect to $\tau$. Let $F_{\tau} \assign \tmop{Id} -K_{\tau}$. A
  solution $( u_{\tau} ,v_{\tau} ,m_{\tau} ,c_{\tau} )$ to system (\ref{eq:
  system with tau}) is a fixed point of $F_{\tau}$, i.e. $F_{\tau} ( u_{\tau}
  ,v_{\tau} ,m_{\tau} ,c_{\tau} ) =0$.

  \subsubsection*{Existence of degree of $F_{\tau}$.}A solution to (\ref{eq: system
  with tau}) for $0 \leqslant \tau \leqslant 1$ exists as soon as the degree
  $F_{\tau}$ is well defined and non zero, \ Let $\Omega \subset X_{a}$ the
  open set
  \[ \Omega = \left\{ ( u,v,m,c ) ;  \| u
     \|_{\mathcal{C}^{1} ( \bar{I}_{a} )} \leqslant M \nocomma ,  \nocomma \|
     v \|_{\mathcal{C}^{1} ( \bar{I}_{a} )} \leqslant M \nocomma ,  \| m
     \|_{\mathcal{C} ( \bar{I}_{a} )} \leqslant M \nocomma ,  \underline{c}
     <c< \bar{c} \right\} \]
  for some positive constants $M$, $\underline{c}$, $\bar{c}$. These constants
  can be chosen so that for all $0 \leqslant \tau \leqslant 1$, $F_{\tau} \neq
  0$ on $\partial \Omega$. Indeed, Proposition \ref{P: a priori estimates}
  provides estimates for the case $\tau =1$, namely $0<c_{\star}
  /2<c<2c^{\star}$, $\| u \|_{\mathcal{C}^{1} ( \bar{I}_{a} )} \leqslant 1+c
  \leqslant 1+c^{\star}$, $\| v \|_{\mathcal{C}^{1} ( \bar{I}_{a} )} \leqslant
  1+ \beta_{2} ( \Lambda ) c^{\star} ( 1+n_{0} m_{u} )$, $\| m \|_{\mathcal{C}
  ( \bar{I}_{a} )} \leqslant m_{u}$. Setting $M \assign \max ( 1+c^{\star} ,1+
  \beta_{2} ( \Lambda )  c^{\star} ( 1+n_{0} m_{u} ) ,m_{u} )$, $\underline{c}
  \assign c_{\star} /2$, $\bar{c} \assign 2c^{\star}$ ensures $F_{1} ( u,v,m,c
  ) \neq 0$ for any $( u,v,m,c ) \subset \partial \Omega$. Notice now that
  rescaling $f_{\tau} \assign \tau f$ and $\phi_{\tau} = \tau \phi$ in
  (\ref{eq: system with tau}), leads to a fixed point problem similar to
  $F_{1}$ but with rescaled reaction and vaporisation terms. The
  conclusions of Proposition \ref{P: a priori estimates} still hold with the
  same bounds, ensuring that for all $0 \leqslant \tau \leqslant 1$,
  $F_{\tau} ( u,v,m,c ) \neq 0$ for any $( u,v,m,c ) \subset \partial \Omega$.

  \subsubsection*{Calculation of the degree of $F_{\tau}$.}Thanks to the properties
  of $K_{\tau}$, and the homotopy invariance of the degree we have $\deg (
  F_{\tau} , \Omega ,0 ) = \deg ( F_{0} , \Omega ,0 )$ for all $0 \leqslant
  \tau \leqslant 1$ and it suffices to compute the degree of $F_{0}$. Given $(
  u,v,m,c )$, the solution $( \hat{u} , \hat{v} , \hat{m} , \hat{c} )$ to
  (\ref{eq: system for the fixed point problem}) with $\tau =0$ is
  \[ \begin{array}{l}
  	\hat{u} ( c ) =e^{c ( x-a )} , \hspace{1em} \hat{v} ( c ) =v_{u} ( 1-e^{c
     	( x-a ) / \Lambda} ) , \medskip\\
	\hat{m} ( u,m,c ) = \mathcal{M} (
     	u,m,c ) , \hspace{1em} \hat{c} ( u,c ) =c-u ( 0 ) + \theta_{i} .
  \end{array} \]
  The mapping
  \[ F_{0} :  ( u,v,m,c ) \rightarrow ( u- \hat{u} ( c ) ,v- \hat{v} ( c ) ,m-
     \hat{m} ( u,m,c ) ,u ( 0 ) - \theta_{i} ) \]
  is homotopic to
  \[ \tilde{F}_{0} :  ( u,v,m,c ) \rightarrow ( u- \hat{u} ( c ) ,v- \hat{v} (
     c ) ,m- \hat{m} ( u,m,c ) , \hat{u} ( 0 ) - \theta_{i} ). \]
  This linear system has a unique solution in $X_{a}$ therefore its degree is
  different from $0$. The same result holds for $F_{1}$ which
  therefore admits at least one solution in $X_{a}$.
\end{proof*}

\subsection{Existence of travelling waves}\label{SS: Existence of travelling
waves}

We are now able to prove the first main result of the paper, that is the
existence of a travelling wave for system (\ref{eq: nondimensional travelling
wave}) on the real line.

\begin{proof*}{Proof of Theorem \ref{T: existence of travelling waves}}
  Let $( u_{a} ,v_{a} ,m_{a} ,c_{a} )$ the solution to (\ref{eq: system on a
  bounded domain}) for $a \geqslant a_{\star}$. We proved that $( u_{a} ,v_{a}
  ,m_{a} ,c_{a} )$ is bounded in $W^{1, \infty} ( -a,a ) \times W^{1, \infty}
  ( -a,a ) \times L^{\infty} ( -a,a ) \times \mathbbm{R}$ uniformly with
  respect to $a \geqslant a_{\star}$. Using system (\ref{eq: system on a
  bounded domain}), the uniform boundedness also holds in $W^{2, \infty} (
  -a,a ) \times W^{2, \infty} ( -a,a ) \times L^{\infty} ( -a,a ) \times
  \mathbbm{R}$. Let $( a_{n} )_{n \geqslant 0}$, with $a_{n} \geqslant
  a_{\star}$ an increasing sequence tending to infinity. We can extract a
  sequence $( u_{a_{n}} ,v_{a_{n}} ,m_{a_{n}} ,c_{a_{n}} )$ converging in
  $\mathcal{C}^{1}_{\tmop{loc}} ( \mathbbm{R} ) \times
  \mathcal{C}^{1}_{\tmop{loc}} ( \mathbbm{R} ) \times
  \mathcal{C}^{0}_{\tmop{loc}} ( \mathbbm{R} ) \times \mathbbm{R}$ solution to
 
\[  \begin{array}{rcl}
    -u'' +c u' & = & f ( u ) v \hspace{1em} \tmop{on}  \mathbbm{R},\\
    - \Lambda v'' +c v' & = & f ( u ) v-c n_{0} m'  _{} \hspace{1em} \tmop{on}
    \mathbbm{R},\\
    c m' & = & - \phi ( u,m ) \hspace{1em} \tmop{on}  \mathbbm{R},\\
    u ( 0 ) & = & \theta_{i}.
  \end{array}\]
 
  From the explicit expressions of $u_{a_{n}}$ and $v_{a_{n}}$ on $( -a_{n} ,0
  )$ it follows that $u ( - \infty ) =0$ and $v ( - \infty ) =v_{u}$.
  Moreover, since the velocity $c_{a_{n}}$ is uniformly bounded below by a
  certain $\underline{c} >0$, the vaporisation point $x_{v}^{a_{n}}$ such that
  $u ( x_{v}^{a_{n}} ) = \theta_{v}$ is bounded below by $x_{v}^{\star} =-
  \log ( \theta_{i} / \theta_{v} ) / \underline{c}$ and $m_{a_{n}} =m_{u}$ on
  $( -a_{n} ,x_{v}^{\star} )$. It follows $m ( - \infty ) =m_{u}$. Now, $v$
  and $v'$ inherit the uniform bounds of $v_{a_{n}}$ and $v_{a_{n}}'$, so that
  $v' ( + \infty )$ is bounded, therefore $v' ( + \infty ) =0$ since $v$ is
  bounded. From the equation satisfied by $v$, since $\phi ( u,m )$ and $f ( u
  ) v$ are bounded, we have also that $v''$ is bounded, therefore $v'' ( +
  \infty ) =0$. Let $x_{v f}^{a_{n}}$ the vaporisation front, that is the
  first point where $m_{a_{n}} ( x_{v f}^{a_{n}} ) =0$. Since $c_{a_{n}}$ is
  bounded above by some $\bar{c} >0$, and since by hypothesis the droplets
  vaporise in finite time, one can bound uniformly $x_{v f}^{a_{n}} <x_{v
  f}^{\star}$ for some $x_{v f}^{\star} >0$, and for all $a_{n}$, $m_{a_{n}} (
  x ) =0$ if $x \geqslant x_{v f}^{\star}$. It follows $m ( \infty ) =0$, and
  the equation for $v$ implies $f ( u ( + \infty ) )  v ( + \infty ) =0$. The
  monotonicity of $u$ holds for the same reason as the monotinicity of
  $u_{a_{n}}$ holds, therefore $u ( + \infty ) > \theta_{i}$ and $f ( u ( +
  \infty ) ) >0$, hence $v ( + \infty ) =0$. Finally, since $m ( x ) =0$ on $(
  x_{v f} ,+ \infty )$, we have from Proposition \ref{P: a priori estimates}
  that $v \geqslant \beta_{1} ( \Lambda ) ( 1-u ) \geqslant 0$ on $( x_{v f}
  ,+ \infty )$ hence $u ( + \infty ) =1$.
\end{proof*}

\section{High activation energy limit}\label{S: High activation energy limit}

\subsection{Characterisation of the limiting profiles}\label{SS:
characterisation of the limiting profiles}

In this section we assume a reaction term $f_{\varepsilon} ( u )$ obeying the
Arrhenius law (\ref{eq: Arrehnius law}) on $[ \theta_{i} ,1 ]$. Set
$\theta_{\varepsilon} \assign 1+A \varepsilon \ln \varepsilon <1$ with $A=100$
and $\varepsilon$ small enough. Notice $\theta_{i} < \theta_{\varepsilon}$ for
$\varepsilon$ small enough. We have immediately
\begin{equation}
  \lim_{\varepsilon \rightarrow 0} \theta_{\varepsilon} =1 \hspace{2em}
  \tmop{and} \hspace{2em} \lim_{\varepsilon \rightarrow 0} \hspace{1em}
  \max_{\theta_{i} \leqslant s \leqslant \theta_{\varepsilon}} f_{\varepsilon}
  ( s ) ( 1-s ) =0, \label{eq: bounds of feps before the reaction zone}
\end{equation}
\begin{equation}
  \lim_{\varepsilon \rightarrow 0} \int_{\theta_{i}}^{1} f_{\varepsilon} ( s )
  ( 1-s ) \mathd s= \lim_{\varepsilon \rightarrow 0} G_{\varepsilon} ( 1 ) =:
  \mu < \infty . \label{eq: total reaction amount}
\end{equation}
We follow and adapt the high energy activation analyses and tools from
{\cite{BerNicSch1985,GlaRoq1994,GlaRoq1996a}}.

\begin{proof*}{Proof of Theorem \ref{T: limiting system in the HEA limit}}
  {\dueto{limiting system in the HAE limit}}Choose $\varepsilon >0$ small
  enough so that $\mu /2 \leqslant G_{\varepsilon} ( 1 ) \leqslant 2 \mu$. The
  bounds (\ref{eq: bounds on the velocity}) for the velocity with $a=+ \infty$
  and $u_{\varepsilon}' ( + \infty ) =0$ allow us to estimate independently of
  $\varepsilon$
  \[ \underline{c} / \sqrt{2} \leqslant c_{\varepsilon} \leqslant \sqrt{2}
     \bar{c} , \]
  with $\underline{c} \assign \min ( c_{\star} ,c_{1} )$, $c_{\star}$ as in
  Proposition \ref{P: bounds on the velocity}, \ $c_{1} = \sqrt{2 \beta_{1} (
  \Lambda ) \mu} / \theta_{i}$, and $\bar{c} \assign \sqrt{2 \beta_{2} (
  \Lambda ) \mu} / \theta_{i}$.

  We can therefore find a decreasing sequence $( \varepsilon_{n} )_{n \in
  \mathbbm{N}}$ such that $c_{\varepsilon_{n}}$ converges
  \[ c_{\varepsilon_{n}} \rightarrow c>0. \]
  For $\varepsilon$ small enough, let $x_{\varepsilon} >0$ be the unique point
  satisfying $u_{\varepsilon} ( x_{\varepsilon} ) =1+A \varepsilon \ln
  \varepsilon =: \theta_{\varepsilon}$<1, where $A$ is a fixed sufficiently
  large constant (take $A=100$). Let $x_{v}^{\varepsilon}$ be the point where
  vaporisation starts, i.e. $u_{\varepsilon} ( x_{v}^{\varepsilon} ) \assign
  \theta_{v}$. Let $x_{v}^{\varepsilon} \leqslant x \leqslant
  x_{\varepsilon}$. Rewriting equation $( \star \star \star )$ on $[
  x_{v}^{\varepsilon} ,x ]$ and comparing $v_{\varepsilon}$ with
  $u_{\varepsilon}$ yields
  \begin{eqnarray*}
    \frac{1}{2} u'_{\varepsilon} ( x )^{2} & \geqslant & \frac{1}{2} 
    c_{\varepsilon}^{2} \theta_{v}^{2} - \int_{0}^{x} f_{\varepsilon} (
    u_{\varepsilon} ) v_{\varepsilon} \cdot u'_{\varepsilon}\\
    & \geqslant & \frac{1}{2}   \underline{c}^{2}   \theta_{v}^{2} - \max (
    1, \Lambda ) \int_{0}^{x} f_{\varepsilon} ( u_{\varepsilon} ) (
    1-u_{\varepsilon} ) \cdot u'_{\varepsilon},
  \end{eqnarray*}
  where the last integral is $\cal{O} ( \varepsilon )$ thanks to the properties of
  $f_{\varepsilon}$. Hence there exist constants $\alpha >0$ and
  $\varepsilon_{0} >0$ such that for any $\varepsilon \leqslant
  \varepsilon_{0}$ and $x \in [ x_{v}^{\varepsilon} ,x_{\varepsilon} ]$, we
  have $u' ( x ) \geqslant \alpha$. Since $u_{\varepsilon} ( 0 ) = \theta_{i}
  < \theta_{\varepsilon} =u_{\varepsilon} ( x_{\varepsilon} )$, it follows
  $0<x_{\varepsilon} \leqslant ( \theta_{\varepsilon} - \theta_{v} ) / \alpha
  < ( 1- \theta_{v} ) / \alpha =:x_{0} <+ \infty$. The uniform boundedness of
  $x_{\varepsilon}$ implies, up to extracting a subsequence,
  \[ x_{\varepsilon_{n}} \rightarrow \bar{x} \leqslant x_{0} . \]
  Now for the point $x_{v}^{\varepsilon}$ where vaporisation starts, we have
  $| x_{v}^{\varepsilon} -0 | \leqslant ( \theta_{v} - \theta_{i} ) / \alpha$,
  hence $x_{v}^{\varepsilon}$ is bounded and converges, after extraction of a
  subsequence,
  \[ x_{v}^{\varepsilon_{n}} \rightarrow x_{v} . \]

  \subsubsection*{Convergence of $u_{\varepsilon}$.}The convergence of
  $c_{\varepsilon}$ together with the explicit expressions of the exponential
  profile $u_{\varepsilon}$ on $( - \infty ,0 )$ implies the convergence in
  $H^{1} ( - \infty ,0 )$ of $u_{\varepsilon}$ towards $u$. Inspecting the
  $L^{2}$ estimates in the proof of Proposition \ref{P: bounds on the
  velocity} yields the uniform $H^{1} ( 0,x_{0} )$ boundedness of
  $u_{\varepsilon}$, therefore convergence in $\mathcal{C}^{0} ( (
  0,x_{0} ) )$. Thanks to the monotonicity of $u_{\varepsilon}$, we have that
  $u_{\varepsilon} \rightarrow 1$ uniformly on $( \bar{x} ,+ \infty )$.
  Finally, the reaction term is uniformly bounded on $( - \infty
  ,x_{\varepsilon} )$ by a multiple (\ref{eq: bounds of feps before the
  reaction zone}), which tends to $0$ by hypothesis. This implies that the
  reaction term tends uniformly towards zero on any compact subset of $( -
  \infty , \bar{x} )$. To summarize, $u_{\varepsilon}$ converges in $H^{1} (
  \mathbbm{R} )$ towards a continuous profile $u$ solution to $-u'' +c u' =0$
  on $( - \infty , \bar{x} )$, $u=1$ on $( \bar{x} ,+ \infty )$.

  \subsubsection*{Convergence of $m_{\varepsilon}$.}Since $\theta_{v} <1$, we have
  for $\varepsilon$ small enough that $x^{\varepsilon}_{v} <x_{\varepsilon}$,
  and $\theta_{\varepsilon} \geqslant ( 1+ \theta_{v} ) /2> \theta_{v}$. From
  the hypothesis on the vaporisation law, it follows that complete
  vaporisation occurs on a finite intervall. Let $x_{{vf}}^{\varepsilon}$
  denote the unique position of the vaporisation front. The uniform
  boundedness of $x_{{vf}}^{\varepsilon}$ implies the convergence
  \[ x_{{vf}}^{\varepsilon_{n}} \rightarrow x_{{vf}}. \]
  Moreover, from the $\mathcal{C}^{0}$ convergence of $u_{\varepsilon}$
  towards $u$ on the compact $[ x_{v} -1,x_{{vf}} +1 ]$, together with
  the Lipshitz properties of the vaporisation terms, it follows that
  $m_{\varepsilon}$ converges towards $m$ on $[ x_{v} -1,x_{{vf}} +1 ]$
  in $\mathcal{C}^{0} ( [ x_{v} -1,x_{{vf}} +1 ] )$ solution to $c m' =-g
  ( u,m )$. Since $m_{\varepsilon}' =0$ outside of $[ x_{v}^{\varepsilon}
  ,x_{{vf}}  ^{\varepsilon} ]$, the convergence holds on $\mathcal{C}^{0}
  ( \mathbbm{R} )$ with $m$ solving $m ( - \infty ) =m_{u}$ and $c m' =-g (
  u,m )$ on $\mathbbm{R}$.

  \subsubsection*{Convergence of $v_{\varepsilon}$.}The convergence of
  $v_{\varepsilon}$ follows the same line as that of $u_{\varepsilon}$ and
  uses the convergence of $m_{\varepsilon}$. The limiting solution $v \in
  \mathcal{C}^{0} ( \mathbbm{R} )$ satisfies $- \Lambda v'' +c v' =-c n_{0}
  m'$on $( - \infty , \bar{x} )$, $v=0$ on $( \bar{x} ,+ \infty )$.

  \subsubsection*{Convergence on $( \bar{x} ,+ \infty )$.}Since $u_{\varepsilon}$ is
  increasing, we have for any $x_{\varepsilon} \leqslant x$, $u (
  x_{\varepsilon_{n}} ) = \theta_{\varepsilon_{n}} \leqslant
  u_{\varepsilon_{n}} ( x ) \leqslant 1$ with $\theta_{\varepsilon_{n}}
  \rightarrow 1$. This implies easily $u_{\varepsilon} ( x ) \rightarrow 1$ on
  $( \bar{x} ,+ \infty )$. The upper bound $v_{\varepsilon} ( x ) \leqslant
  \beta_{2} ( \Lambda ) ( 1-u_{\varepsilon} ( x ) )$ ensures $v_{\varepsilon}
  ( x ) \rightarrow 0$ on $( \bar{x} ,+ \infty )$. Now, the auxiliary function
  $y_{\varepsilon} \assign u_{\varepsilon} + \Lambda v_{\varepsilon} -1$
  satisfies $y_{\varepsilon} ( x_{\varepsilon} ) = \theta_{\varepsilon} +
  \Lambda v_{\varepsilon} ( x_{\varepsilon} ) -1$, $y_{\varepsilon} ( + \infty
  ) =0$, and $-y_{\varepsilon}'' +c_{\varepsilon} y_{\varepsilon}' = ( \Lambda
  -1 ) c_{\varepsilon} v'_{\varepsilon} - \Lambda c_{\varepsilon} n_{0}
  m'_{\varepsilon}$. Therefore the $H^{1}$--limit $y$ satisfies $y=u+ \Lambda
  v-1=0$ on $( \bar{x} ,+ \infty )$, together with $-y'' +c y' = ( \Lambda -1
  ) c v' - \Lambda c n_{0} m'$. This implies $m' =0$ on $( \bar{x} ,+ \infty )
  \nocomma$ therefore also $m ( x ) =0$ on $( \bar{x} ,+ \infty )$. Obviously,
  $u$ and $v$ are $\mathcal{C}^{2}$ on $( \bar{x} ,+ \infty )$.

  \subsubsection*{Estimate of $x_{{vf}}$.}The fact that $m ( x ) =0$ for $x
  \geqslant \bar{x}$ is equivalent to
  \[ x_{{vf}} \leqslant \bar{x}, \]
  that is in the HAE limit, the vaporisation region ends before or at the
  reaction zone.

  \subsubsection*{Convergence on $( - \infty , \bar{x} )$.}Since $u_{\varepsilon}$
  is increasing and $u_{\varepsilon} ( 0 ) = \theta_{i}$, the reaction term is
  zero on $( - \infty ,0 )$. From property (\ref{eq: bounds of feps before the
  reaction zone}) together with the upper estimate $v_{\varepsilon} \leqslant
  \beta_{2} ( \Lambda ) ( 1-u_{\varepsilon} )$, we deduce that
  $f_{\varepsilon} ( u_{\varepsilon} ) v_{\varepsilon}$ converges uniformly to
  zero on any compact set of $[ 0, \bar{x} )$, therefore also on any compact
  set of $( - \infty , \bar{x} )$. This implies that the limiting functions
  $u$, $v$ and $m$ satisfiy $-u'' +c u' =0$, $- \Lambda v'' +c v' =-c n_{0}
  m'$, and $c n_{0} m' =-g ( u,m )$ on $( - \infty , \bar{x} )$, with $m' \in
  L^{2} ( ( - \infty , \bar{x} ) )$, $v \in H^{2} ( ( - \infty , \bar{x} ) )
  \subset \mathcal{C}^{1} ( ( - \infty , \bar{x} ) )$ and $u ( x ) =
  \theta_{i} e^{c ( x- \bar{x} )} \in \mathcal{C}^{2} ( ( - \infty , \bar{x} )
  )$.

  This concludes the proof of Theorem \ref{T: limiting system in the HEA
  limit} {\dueto{limiting system in the HAE limit}}
\end{proof*}

\begin{proof*}{Proof of Theorem \ref{T: rigorous internal layer analysis}}
  {\dueto{internal layer analysis}}Set now $u_{\varepsilon} ( 0 ) =1+A
  \varepsilon \ln \varepsilon =: \theta_{\varepsilon}$ at $x=0$. Let then
  $x_{i}^{\varepsilon} \leqslant 0$ the ignition point where $u_{\varepsilon}
  ( x_{i}^{\varepsilon} ) = \theta_{i}$. We have seen in the proof of Theorem
  \ref{T: limiting system in the HEA limit} above that $| x_{i}^{\varepsilon}
  -0 |$ is bounded independently of $\varepsilon$. From property (\ref{eq:
  bounds of feps before the reaction zone}), we have thereofre
  $f_{\varepsilon} ( u_{\varepsilon} ) v_{\varepsilon} = \mathcal{O} (
  \varepsilon^{A/2} )$ uniformly on $( x_{i}^{\varepsilon} ,0 )$ and is zero
  on $( - \infty ,x_{i}^{\varepsilon} )$. Integrating the equation for
  $u_{\varepsilon}$ between $- \infty  $ and $0$, this implies
  \[ -u'_{\varepsilon} ( 0 ) +c_{\varepsilon} \theta_{  \varepsilon} =
     \mathcal{O} ( \varepsilon^{A/2} ) , \]
  where $\theta_{\varepsilon} =1+ \mathcal{O} ( \varepsilon \ln \varepsilon
  )$. It follows that
  \[ c_{\varepsilon} =u'_{\varepsilon} ( 0 ) + \mathcal{O} ( \varepsilon \ln
     \varepsilon ) . \]
  \subsubsection*{Rescaled system.}Let $\xi \assign x/ \varepsilon$.
  \[ \hat{u}_{\varepsilon} ( \xi ) = \frac{u_{\varepsilon} ( x )
     -1}{\varepsilon} , \hspace{1em} \hat{v}_{\varepsilon} ( \xi ) =
     \frac{v_{\varepsilon} ( x )}{\varepsilon} , \hspace{1em} m_{\varepsilon}
     ( \xi ) =\hat{m}_{\varepsilon} ( x ) , \]
  together with the usual auxiliary functions
  \[ z_{\varepsilon} ( x ) =u_{\varepsilon} ( x ) +v_{\varepsilon} ( x ) -1,
     \hspace{1em} y_{\varepsilon} ( x ) =u_{\varepsilon} ( x ) + \Lambda
     v_{\varepsilon} ( x ) -1, \]
  \[ \hat{z}_{\varepsilon} ( \xi ) = \frac{z_{\varepsilon} ( x )}{\varepsilon}
     , \hspace{1em} \hat{y}_{\varepsilon} ( \xi ) = \frac{y_{\varepsilon} ( x
     )}{\varepsilon} . \]
  These functions obey the system of equations:
  \[ - \hat{u}_{\varepsilon}'' + \varepsilon c_{\varepsilon}
     \hat{u}_{\varepsilon}' = \hat{v}_{\varepsilon}   \exp (
     \hat{u}_{\varepsilon} ) = ( \hat{z}_{\varepsilon} - \hat{u}_{\varepsilon}
     ) \exp ( \hat{u}_{\varepsilon} ) = \frac{1}{\Lambda} (
     \hat{y}_{\varepsilon} - \hat{u}_{\varepsilon} ) \exp (
     \hat{u}_{\varepsilon} ) , \]
  \[ \hspace{1em} \hat{u}_{\varepsilon} ( 0 ) =A  \ln \varepsilon ,
     \hspace{1em} \hat{u}_{\varepsilon} ( + \infty ) =0,  \hspace{1em}
     \hat{u}_{\varepsilon} ( - \infty ) =-1/ \varepsilon , \]
  \[ - \hat{v}_{\varepsilon}'' + \varepsilon c_{\varepsilon}
     \hat{v}_{\varepsilon}' =- \exp ( \hat{u}_{\varepsilon} )  
     \hat{v}_{\varepsilon} -c_{\varepsilon} n_{0} \hat{m}' , \hspace{1em}
     \hat{v}_{\varepsilon} ( + \infty ) =0, \hspace{1em} \hat{v}_{\varepsilon}
     ( - \infty ) =v_{u} / \varepsilon, \]
  \[ c_{\varepsilon} n_{0} \hat{m}_{\varepsilon}' =-g ( 1+ \varepsilon
     \hat{u}_{\varepsilon} ) \phi ( \hat{m}_{\varepsilon} ) , \hspace{1em}
     \hat{m}_{\varepsilon} ( - \infty ) =m_{u} , \]
  \[ - \hat{z}_{\varepsilon}'' + \varepsilon c_{\varepsilon}
     \hat{z}_{\varepsilon}' = ( \Lambda -1 ) \hat{v}_{\varepsilon}''
     -c_{\varepsilon} n_{0} \hat{m}_{\varepsilon} ,
   \]
   \[  \hspace{5cm} \hat{z}_{\varepsilon} ( + \infty ) =0, \hspace{1em} \hat{z}_{\varepsilon}
     ( - \infty ) = ( v_{u} -1 ) / \varepsilon, \]
  \[ - \Lambda \hat{y}_{\varepsilon}'' + \varepsilon c_{\varepsilon}
     \hat{y}_{\varepsilon}' = \varepsilon c_{\varepsilon} ( 1- \Lambda )
     \hat{u}_{\varepsilon}' -c_{\varepsilon} n_{0} \hat{m}_{\varepsilon} ,
 \]
 \[    \hspace{5cm} \hat{y}_{\varepsilon} ( + \infty ) =0, \hspace{1em}
     \hat{y}_{\varepsilon} ( - \infty ) = ( \Lambda v_{u} -1 ) / \varepsilon.
  \]
  Notice that we have $u_{\varepsilon}' ( 0 ) = \hat{u}_{\varepsilon}' ( 0 )$,
  therefore also $c_{\varepsilon} = \hat{u}'_{\varepsilon} ( 0 ) + \mathcal{O}
  ( \varepsilon \ln \varepsilon )$. Our goal is to estimate
  $\hat{u}'_{\varepsilon} ( 0 )$ in the limit $\varepsilon \rightarrow 0$.

  \subsubsection*{Approximate system.}Let $\tilde{u}_{\varepsilon}$ and
  $\tilde{y}_{\varepsilon}$ the solution to the linear system

  \begin{equation}
    - \tilde{u}_{\varepsilon}'' = \frac{1}{\Lambda} ( \tilde{y}_{\varepsilon}
    - \tilde{u}_{\varepsilon} ) \exp ( \hat{u}_{\varepsilon} ) , \label{eq:
    approximate system u}
  \end{equation}
  \begin{equation}
    - \Lambda \tilde{y}_{\varepsilon}'' =-c_{\varepsilon} n_{0}
    \hat{m}'_{\varepsilon}, \label{eq: approximate system y}
  \end{equation}
  with the same boundary conditions as $\hat{u}_{\varepsilon}$ and
  $\hat{y}_{\varepsilon}$, and where we omitted the terms involving the
  factors $\varepsilon c_{\varepsilon}$. Classical results in elliptic
  regularity ensure that $( \hat{u}_{\varepsilon} , \hat{y}_{\varepsilon} ) -
  ( \tilde{u}_{\varepsilon} , \tilde{y}_{\varepsilon} )$ is of order
  $\mathcal{O} ( \varepsilon )$ on $[ 0, \infty )$ in $\mathcal{C}^{1}$.
  Therefore $c_{\varepsilon} = \tilde{u}'_{\varepsilon} ( 0 ) + \mathcal{O} (
  \varepsilon \ln \varepsilon )$.
  
  \subsubsection*{Case $x_{{vf}}<0$.}Let $x_{{vf}}^{\varepsilon}$ the
  position of the vaporisation front such that $\hat{m}_{\varepsilon} =0$ on $(
  x_{{vf}}^{\varepsilon} ,+ \infty )$. If $x_{{vf}} <0$, then for
  $\varepsilon$ small enough, $\hat{m}_{\varepsilon} =0$ and
  $\tilde{y}_{\varepsilon} =0$ on $( 0,+ \infty )$ and
  
  \[ \begin{array}{lll}
  - \tilde{u}_{\varepsilon}'' &=& - \tilde{u}_{\varepsilon} \exp (
     \hat{u}_{\varepsilon} ) / \Lambda =- \tilde{u}_{\varepsilon} \exp (
     \tilde{u}_{\varepsilon} ) / \Lambda + \tilde{u}_{\varepsilon}  
     \mathcal{O ( \varepsilon )} 
     \medskip\\
     &=& - \tilde{u}_{\varepsilon} \exp (
     \tilde{u}_{\varepsilon} ) / \Lambda + \mathcal{O ( \varepsilon \ln
     \varepsilon )} . 
     \end{array}
     \]
  Omitting the last term in the equation above, it remains to study
  the solution to the new approximate system with the same boundary conditions
  \[ - \check{u}_{\varepsilon}'' =- \check{u}_{\varepsilon} \exp (
     \check{u}_{\varepsilon} ) / \Lambda. \]
  Integrating against $\check{u}_{\varepsilon}'$ on $( 0,+ \infty )$ yields
  \[ \frac{\check{u}_{\varepsilon}' ( 0 )^{2}}{2} = \frac{\mu}{\Lambda}
     \int_{A  \ln \varepsilon}^{0} - \sigma e^{\sigma}, \]
  where the integral in the right hand side tends to one as $\varepsilon$ goes to zero. Hence
  \[ c_{\varepsilon} =u_{\varepsilon}' ( 0 ) = \check{u}_{\varepsilon}' ( 0 )
     + \mathcal{O} ( \varepsilon \ln \varepsilon ) = \sqrt{\frac{2
     \mu}{\Lambda}} + \mathcal{O} ( \varepsilon \ln \varepsilon ) , \]
  and the limiting velocity is $c= \sqrt{2 \mu / \Lambda}$ when $x_{{vf}}
  =0$.
  
  \subsubsection*{Case $x_{{vf}} =0$.}In the case where the vaporisation front
  is located in the limit at the position $\bar{x} =0$ of the reaction zone,
  the velocity $c$ is imposed by the condition $m ( x ) =0$. Indeed, we have
  in the limit $u ( x ) = e^{c x}$ for $x \leqslant 0$. Assume $0<c_{1}
  <c_{2}$. The corresponding profiles satisfy $u_{1} >u_{2}$ on $( - \infty ,0
  )$ and the corresponding positions where vaporisation starts satisfy
  $x_{v}^{1} <x_{v}^{2}$. The monotonicity properties of the vaporisation law
  imply $x_{{vf}}^{1} <x_{{vf}}^{2}$, where those are the
  vaporisation fronts associated to the velocities $c_{1}$ and $c_{2}$
  respectively. Consequently, the vaporisation front $x_{{vf}}$ is an
  increasing function $x_{{vf}} \assign x_{{vf}} ( c )$ of the
  velocity $c$. With similar arguments, it is also an increasing function of
  $m_{u}$. Moreover, since the vaporisation law $\phi$ is Lipschitz on $(\theta_v,1)$, it follows that $x_{{vf}}
  ( c )$ is continuous and we can define in a unique manner $c_{\star} ( m_{u} )$
  to be the velocity such that $x_{{vf}} ( c_{\star} ( m_{u} ) ) =0$.

  As a conclusion, when $x_{{vf}} =0$, we have $c=c_{\star} ( m_{u} )$,
  whereas $c= \sqrt{2 \mu / \Lambda}$ for $x_{{vf}} <0$. The monotonicity
  of $x_{{vf}} :c \mapsto x_{{vf}} ( c )$ implies
  \[ c_{\star} ( m_{u} ) \leqslant \sqrt{2 \mu / \Lambda}, \]
  \[ x_{{vf}} =0 \hspace{1em} \Leftrightarrow \hspace{1em} c=c_{\star} (
     m_{u} ) < \sqrt{2 \mu / \Lambda}, \]
  \[ x_{{vf}} <0 \hspace{1em} \Leftrightarrow \hspace{1em} c= \sqrt{2 \mu
     / \Lambda}. \]
  This concludes the proof of Theorem \ref{T: rigorous internal layer
  analysis}.
\end{proof*}

\subsection{Internal combustion layer in the vaporisation controlled
regime}\label{SS: internal layer in the vaporisation limited regime}

We are interested in the case where the droplets finish vaporising inside the
combustion layer. In order to measure the quantity of droplets still present
in the combustion region $[0,+\infty)$, it is convenient to introduce $M_{\varepsilon} ( x
)$ the following primitive of the mass of liquid
\[ M_{\varepsilon} ( x ) \assign \int_{x}^{+ \infty} m_{\varepsilon} ( x )  
   \mathd x. \]

\begin{theorem}
  \label{T: estimate of the overlapping region}{\dueto{estimate of the
  overlapping region}}Let $M_{\varepsilon}$ defined above. We have
  \[ M_{\varepsilon} ( 0 ) = \mathcal{O} ( \varepsilon   \ln \varepsilon ) .
  \]
\end{theorem}

\begin{proof*}{Proof.}
  Consider the approximate system on the intervall $( 0,+ \infty )$

  \begin{equation}
    - \check{u}_{\varepsilon}'' = \frac{1}{\Lambda} ( \check{y}_{\varepsilon}
    - \check{u}_{\varepsilon} ) \exp ( \check{u}_{\varepsilon} ) ,
    \hspace{1em} \check{u} ( 0 ) =A  \ln \varepsilon , \hspace{1em}
    \check{u}_{\varepsilon} ( + \infty ) =0, \label{eq: approximate system u
    overlapping region}
  \end{equation}
  \begin{equation}
    - \Lambda \check{y}_{\varepsilon}'' =-c_{\varepsilon} n_{0}
    \hat{m}'_{\varepsilon} , \hspace{1em} \check{y}_{\varepsilon} ( 0 ) =
    \hat{y}_{\varepsilon} ( 0 ) , \hspace{1em} \check{y}_{\varepsilon} ( +
    \infty ) =0, \label{eq: approximate system y overlapping region}
  \end{equation}
  such that $\check{u}_{\varepsilon} - \hat{u}_{\varepsilon} = \mathcal{O} (
  \varepsilon )$ and $\check{y}_{\varepsilon} - \hat{y}_{\varepsilon} =
  \mathcal{O} ( \varepsilon )$. We also rescale $M_{\varepsilon} ( x )$,
  \[ \hat{M}_{\varepsilon} ( \xi ) \assign \int_{\xi}^{+ \infty}
     \hat{m}_{\varepsilon} ( \xi )   \mathd \xi . \]
  Let $\xi_{{vf}}^{\varepsilon}$ the position of the vaporisation front
  and assume $\xi_{{vf}}^{\varepsilon} >0$, that is some vaporisation
  occurs after $\xi =0$. We have $\hat{m}'_{\varepsilon} =0$ on $(
  \xi^{\varepsilon}_{{vf}}  ,+ \infty )$, hence $\check{y}_{\varepsilon}$
  is linear with limit zero at infinity, therefore $\check{y}_{\varepsilon}
  =0$ on $( \xi^{\varepsilon}_{{vf}}  ,+ \infty )$. It follows that
  \[ - \Lambda \check{y}_{\varepsilon} ( \xi ) =+c_{\varepsilon} n_{0}
     \hat{M}_{\varepsilon} ( \xi ) . \]
  On the other hand, we had for the original rescaled unknown
  \[ \hat{y}_{\varepsilon} ( 0 ) = \hat{u}_{\varepsilon} ( 0 ) + \Lambda
     \hat{v}_{\varepsilon} ( 0 ) \ln \varepsilon, \]
  because $\hat{u}_{\varepsilon} ( 0 ) =A  \ln \varepsilon$, and (Proposition
  \ref{P: qualitative properties} \& \ref{P: a priori estimates}), $0
  \leqslant \hat{v}_{\varepsilon} \leqslant \beta_{2} ( \Lambda ) |
  \hat{u}_{\varepsilon} |$. It follows
  \[ - \Lambda \check{y}_{\varepsilon} ( 0 ) =- \Lambda \hat{y}_{\varepsilon}
     ( 0 ) + \mathcal{O} ( \varepsilon ) = \mathcal{O} ( \varepsilon
     \mathcal{\ln \varepsilon} ) . \]
  Since $\hat{M}_{\varepsilon} ( 0 ) =M_{\varepsilon} ( 0 ) / \varepsilon$, we
  proved
  \[ M_{\varepsilon} ( 0 ) = \mathcal{O} ( \varepsilon   \ln \varepsilon ) \label{eq: estimate of M0}.
  \]
  This concludes the proof of Theorem \ref{T: estimate of the overlapping
  region}.
  \end{proof*}

\subsubsection*{Application to the $d^{2}$--law.}Assume a vaporisation term of the
form
\[ \forall \theta_{v} \leqslant u \leqslant 1, \hspace{1em} \phi ( u,m )
   =g_{0}  m^{\delta} , \hspace{1em} g_0>0, \hspace{1em} 0 \leqslant \delta <1. \]
Integrating twice on $[ x_{\varepsilon} ,x^{\varepsilon}_{{vf}} ]$ yields
\[ \forall x \geqslant x^{\varepsilon}_{v} , \hspace{1em} m ( x ) =
   \frac{g_{0}}{c n_{0}} \hspace{1em} ( 1- \delta )^{1/ ( 1- \delta )}
   \hspace{1em} ( x^{\varepsilon}_{{vf}} -x )^{1/ ( 1- \delta )}, \]
\[ \forall x \geqslant x_{v} , \hspace{1em} M_{\varepsilon} ( x ) =
   \frac{g_{0}}{c n_{0}} \hspace{1em} \frac{( 1- \delta )^{( 2- \delta ) / (
   1- \delta )}}{2- \delta} \hspace{1em} ( x^{\varepsilon}_{{vf}} -x )^{(
   2- \delta ) / ( 1- \delta )}. \]
The estimate (\ref{eq: estimate of M0}) above implies, as soon as
$x_{{vf}}^{\varepsilon} >0$ that
\[ x_{{vf}}^{\varepsilon} = \mathcal{O} ( \varepsilon \ln \varepsilon
   )^{( 1- \delta ) / ( 2- \delta )}. \]
The so--called ``$d^{2}$--law'' states that the rate of variation of the
surface area of a vaporising droplet is approximately constant and corresponds
to the case $\delta =1/3$. We let the reader check that it satisfies the
hypothesis imposed on the vaporisation term, in particular that complete
vaporisation occurs in finite time. The equation above reads
\[ x_{{vf}}^{\varepsilon} = \mathcal{O} ( \varepsilon \ln \varepsilon
   )^{2/5} . \]
This should be compared with the size of the reaction zone in the case
$x_{{vf}} <0$, where $x_{{vf}}^{\varepsilon} \leqslant \gamma <0$ for
some $\gamma$ and $\varepsilon$ small enough. In that case, the spray
travelling wave has the same velocity as the purely gaseous flame with same
temperature in the burnt gas. Also, the analysis of the internal combustion
layer is similar to the gaseous case and provides an estimate of $\mathcal{O (
\varepsilon \ln \varepsilon )}$ for the size of the reaction zone. This
analysis suggests that the internal combustion layer in the presence of droplets
might be considerably stretched and provides a quantitative upper estimate in the
high activation energy limit. This is in accordance with many observations
in the literature. See in particular the numerical
experiments in {\cite{SuaNicHal2001}} that illustrate the fact that the spray
flame consists of a sharp deflagration front followed by a longer afterburn
region where the droplets continue to be slowly vaporised.

\section{Extensions}\label{S: extensions}

In this section, we state without proofs some easy extensions of our analysis.
We consider (i) the case where vaporisation is instantaneous, (ii) the case where the
spray is polydisperse, and (iii) other geometries. We refer to
{\cite{Ber2003}} for more details and extensions.

\subsubsection*{Heuristic analysis in the high activation energy limit.}
A feature specific to our model is the existence of a simple alternative in the high
activation energy limit, where it is easy to predict whether the system is in the
diffusion or vaporisation controlled regime. Assume the system is in the
diffusion controlled regime. This implies that the velocity of the flame is
equal to that of the equivalent gaseous flame with same temperature in the
burnt gas. The temperature profile in that case is known, hence also the
corresponding liquid phase profile, from which one deduces the position
$x_{{vf}}$ where the vaporisation ends. If the vaporisation ends in the
preheating zone before the reaction front, i.e. $x_{{vf}} \leqslant \bar{x}$,
then we are indeed in the diffusion controlled regime. If on the contrary $x_{{vf}}>\bar{x}$, this is a
contradiction with our analysis and the spray flame is in the vaporisation
controlled regime. In that case the position $\bar{x}$ of the reaction
front obeys the constraint \ $\bar{x} =x_{{vf}}$, which
corresponds to a unique value of $c$.

Notice also that thanks to the monotonicity of the vaporisation law we have
seen that the transition from diffusion to vaporisation controlled regime is
sharp and occurs at a certain $m_{u} =m_{u}^{\star}$. As an experience of
thought, fix $\rho_{l}$, $0< \rho_{l} <1$, and vary $m_{u}$ from zero to
infinity and $n_{0}$ accordingly so that $n_{0} m_{u} = \rho_{l}$ for all
values of $m_{u}$. Set $v_{u} =1- \rho_{l}$, and solve for the system in the
HAE limit. Our analysis shows that $c= \sqrt{2 \mu / \Lambda}$ for $m_{u} \in
( 0,m_{u}^{\star} )$ and $c \assign c_{\star} ( m )$ for $m>m_{u}^{\star}$, a
decreasing function of $m$. This is in accordance with the
numerical experiments shown in {\cite{SuaNicHal2001}}.

\subsubsection*{Fast vaporisation.}The situation where the vaporisation is very
rapid can be modeled by an instantaneous vaporisation (see
{\cite{ClaSun1990}}) leading to a Dirac model for the vaporisation term
\[ \phi ( u,m ) =m_{u}   \delta ( x=x_{v} ). \]
Since $\theta_{v} < \theta_{i}$, vaporisation ends before the reaction zone at
$x_{v f} =x_{v} < \bar{x}$. Explicit expressions of the profiles can be found.

\subsubsection*{Polydisperse sprays.}All the previous results can be extended in the
case of polydisperse sprays, that is the situation where droplets of possibly
different sizes are present at a given position. In that case, we need a
statistical description of the distribution of size of the droplets as pointed
out by \tmname{Williams} {\cite{Wil1958a,Wil1988}}. See also \tmname{O'Rourke} \cite{ORo1981}.
For our problem, we assume that the
distribution of size of droplets is known in the fresh gases.

Let $\mathd n ( t,x,m ) = \nu ( t,x,m ) \mathd m$ denote the density number of
droplets at $( t,x )$ with size in the intervall $[ m,m+ \mathd m ]$. The
total mass density $\rho_{l} ( t,x )$ of the liquid phase is
\[ \rho_{l} ( t,x ) = \int_{m} \nu ( t,x,m ) m  \mathd m. \]
and the function $\nu ( t,x,m )$ obeys the conservation equation
\[ \partial_{t} \nu + \partial_{m} ( - \phi ( u,m )   \nu ) =0 \nocomma . \]
The total vaporisation rate of the liquid phase is
\[ \int_{m} \nu ( t,x,m )   \phi ( u,m )   \mathd m, \]
and the corresponding coupled system writes,
\[ \begin{array}{rll}
     -u'' +c u' & = & f ( u ) v \hspace{1em} \tmop{on}  \mathbbm{R},\\
     - \Lambda v'' +c v' & = & -f ( u ) v+  \displaystyle\int_{m} \nu ( t,x,m )  _{} \phi (
     u,m ) \mathd m \hspace{1em} \tmop{on}  \mathbbm{R},\\
     c  \nu' + \partial_{m} ( - \phi ( u,m ) \nu ) & = & 0 \hspace{1em}
     \tmop{on}  \mathbbm{R},\\
     u ( - \infty ) =0, &  & u ( + \infty ) =1,\\
     v ( - \infty ) =v_{u} , &  & v ( + \infty ) =0,\\
     \nu ( - \infty ) = \nu_{u} ( m ) . &  & 
   \end{array} \]
where $\nu_{u} ( m )$ is the distribution profile of mass of droplets in the
fresh gas.

This polydisperse model admits travelling waves (see also
{\cite{LauMas2002a}}). In the HAE limit, the regime of the flame is determined
by the size of the largest droplet in the fresh gas. This result is a
consequence of the simplicity of our vaporisation law, where no collective
effect, nor influence of the density of the gaseous reactant, nor the fine
geometry of the system at small scales were taken into account.

\subsubsection*{Other geometries.}The large activation energy limit analysis can be
carried out in other onedimensionnal geometries, such as the anchored flame on
the half-line, counter-flow-like configurations, or finally radial
geometries in dimension $2$ or $3$. In each case, on can solve explicitly the
problem (provided the expression of the vaporisation law is relatively
simple), and it is possible to study the influence of different parameters on
the existence of a profile or on the value of the burning rate. Relevant
parameters are the typical velocities of injection of the gas and/or droplets,
the value of the vaporisation rate, or the space dimension. We refer the
reader to {\cite{Ber2003}} for more details.

\section{Perspectives and Concluding Remarks}

This paper considered the existence of travelling fronts for a simple
onedimensionnal thermo-diffusive lean spray flame model. We proved the
existence of travelling waves for a class of combustion and vaporisation laws.
As far as qualitative results are concerned, the most important part of the
paper is the study of the high activation energy limit for the system. The
limiting problem involves simple explicit profiles and preserves some
important features of the dynamics. Extensions of these results to the cases
of fast vaporisation, polydisperse sprays or other onedimensionnal geometries
were briefly mentionned.

The present work is a first a step towards rigorous derivation of new
asymptotic models for spray flames. However, much remains to be done in order
to understand the effect of droplets on the dynamics of the flame structure,
as observed by many physicists and experimentalists (see
{\cite{GreMcIBri1999}} and the references therein).
Indeed, our model omits among other things the effects of the latent
heat or more complex couplings in the vaporisation law (see
{\cite{LauMas2002a}} for such considerations). It would be very interesting to derive some high
activation energy models incorporating these effects, in the hope of deriving
explicit expressions of the limiting profiles and combustion rates, as well as
extinction limits. Also, a striking consequence of the high activation energy
asymptotics is that droplets {\tmem{cannot}} cross the flame front in that
limit, but may enter the reaction zone only for large but not infinite values
of the activation energy. It would be interesting to derive intermediate
asymptotic models providing a precise understanding of the structure and
thickening of the combustion region in the presence of vaporising droplets.
Ideally, one would hope to recover some of the features of spray flames
described in
{\cite{Gre2003,NicHalSua2007a,KagGreSiv2009a,KagGreSiv2012a,NicDenHal2014a}}.

As far as dynamical phenomena are concerned, the problem of the stability of
spray flame systems is crucial. The work {\cite{ConDomRoqRyz2006a}} presents
a mathematical analysis of spray pulsating waves. See also {\cite{NicHalSua2005a,Ber2003,BerDomRoq2005}}.
As a second step we have in
mind the problem of acoustic instabilities in spray flame systems. In this
direction, let us mention the work of Clavin an Sun {\cite{ClaSun1990}}. Their
mathematical analysis relies on the possibility of deriving explicit
expressions of the solutions to the problem. This is only possible in certain
asymptotic limits. For that, the authors consider the large activation energy
limit for combustion phenomena and assume an instantaneous vaporisation of the
droplets. As a consequence, both the combustion and the vaporisation zones
reduce to infinitely small regions, whose internal layer is analysed.

However, assuming that the vaporisation zone is very small compared to the
preheating zone is a very restrictive assumption in many applications where
the droplets can spread into the preheating zone, approach the combustion
zone, or enter the combustion zone. We have shown in the present paper that it
is possible to analyse such situations where the droplets approach
or reach the reaction front, and where the slowly vaporising droplets induce a
dramatic change of the combustion rate. We hope that our work will motivate
new studies carrying out rigorous mathematical analysis of asymptotic models
for spray flames.

\section*{Acknowledgments}
It is our pleasure to thank J.-M. Roquejoffre for
suggesting the study and many helpful advices, as well as P.
Haldenwang for very interesting and stimulating discussions.

\end{document}